\documentclass[final]{siamltex}

\usepackage{epsfig,amssymb,latexsym}
\usepackage{amsfonts,psfrag,amsmath,bbm,color}
\usepackage{cancel}

\usepackage{mathrsfs}
\usepackage{graphicx}
\usepackage{textcomp}
\usepackage{multirow}
\usepackage{enumerate}
\usepackage{cancel}
\usepackage{algpseudocode}
\usepackage{caption}  
\usepackage{url}

\usepackage{algorithm}
\usepackage{placeins}

\usepackage{tikz}
\usetikzlibrary{shapes.geometric, arrows}
\tikzstyle{startstop} = [rectangle, rounded corners, minimum width=1cm, minimum height=1cm,text centered, draw=black]
\tikzstyle{io} = [trapezium, trapezium left angle=70, trapezium right angle=110, minimum width=1cm, minimum height=1cm, text centered, draw=black, fill=blue!30]
\tikzstyle{method} = [rectangle, rounded corners, minimum width=1cm, minimum height =1cm, text centered, draw=black]
\tikzstyle{process} = [rectangle, minimum width=1cm, minimum height=1cm, text centered, draw=black]
\tikzstyle{decision} = [diamond, minimum width=0.5cm, minimum height=0.5cm, text centered, draw=black, fill=green!30]
\tikzstyle{arrow} = [thick,->,>=stealth]

\usepackage{amscd}

\graphicspath{{./figs/}}

\newtheorem{remark}{Remark}[section]

\def\PP{{{\rm l}\kern - .15em {\rm P} }}
\def\PN2{{\PP_{N}-\PP_{N-2}}}

\newcommand{\R}{\mathbbm{R}}

\newcommand{\orr}[1]{\overline{#1}^r}

\newcommand{\deleted}[1]{{}}

\begin{document}
\title{Nonintrusive Stabilization of Reduced Order Models for Uncertainty Quantification of Time-Dependent Convection-Dominated Flows}

\author{
	M. Gunzburger
	\thanks{Department of Scientific Computing, Florida State University, \textcolor{black}{Tallahassee, FL 32306, supported by DOE grant DE-SC0009324 and AFSOR grant FA9550-15-1-0001, email: mgunzburger@fsu.edu.}}
	\and T. Iliescu
	\thanks{Department of Mathematics, Virginia Tech, Blacksburg, VA 24061,
		supported by NSF grant \textcolor{black}{ DMS-1821145}, email: iliescu@vt.edu.}
	\and M. Mohebujjaman
	\thanks{Department of Mathematics, Virginia Tech, Blacksburg, VA 24061,
		email: jaman@vt.edu.}
	\and M. Schneier
	\thanks{Department of Mathematics, University of Pittsburgh, Pittsburgh, PA 15260, \textcolor{black}{supported by DOE grant DE-SC0009324 and AFSOR grant FA9550-15-1-0001,} email: mhs64@pitt.edu.}
}

\date{\today}

\maketitle

\begin{abstract}
In this paper, we propose a nonintrusive filter-based stabilization of reduced order models (ROMs) for uncertainty quantification (UQ) of the time-dependent Navier-Stokes equations in convection-dominated regimes.
We propose a novel high-order ROM differential filter and use it in conjunction with an evolve-filter-relax algorithm to attenuate the numerical oscillations of standard ROMs. We also examine how stochastic collocation methods (SCMs) can be combined with the evolve-filter-relax algorithm for efficient UQ of fluid flows.We emphasize that the new stabilized SCM-ROM framework is nonintrusive and can be easily used in conjunction with legacy flow solvers.
We test the new framework in the numerical simulation of a two-dimensional flow past a circular cylinder with a random viscosity that yields a random Reynolds number with mean $Re=100$.
\end{abstract}



\pagestyle{myheadings}
\thispagestyle{plain}
\markboth{\uppercase{M. Gunzburger, T. Iliescu, M. Mohebujjaman, AND M. Schneier}}{\uppercase{Nonintrusive Stabilization of Reduced Order Models}}

\section{Introduction}
   
Mathematical models of realistic flows inherently rely on fundamental input quantities whose actual values are imperfectly measured or sensitive to variations in other parameters. These quantities include the initial conditions, forcing functions, and model coefficients. Robust flow models must identify and properly handle parameter influence whose sensitivity or uncertainty affects model behavior. Such robust flow models are vitally important for numerical weather forecasting, ocean modeling, aeronautical design, and hemodynamics. In this paper, we consider the uncertainty quantification (UQ) of a fluid flow model in the form of the incompressible Navier-Stokes Equations (NSE).

In order to recover accurate solution statistics for fluid flow models, several challenges must be tackled:
(i) These flow models must be realized over an ensemble of parameters at a high spatial resolution. For many flows individual runs can take on the order of hours, which will often make this problem computationally intractable;
(ii) Realistic flows generally occur in the convection-dominated regime, for which standard numerical methods usually yield numerical oscillations; and 
(iii) Legacy codes are often used in realistic fluid flow applications in engineering and geophysics.

A popular approach for dealing with issue (i) is stochastic collocation methods (SCMs) \cite{BNT07,CCMFT2013,gunzburger2014stochastic,FTW2008,stoyanov2013hierarchy,stoyanov2016dynamically}. These methods use global polynomial approximations, taking advantage of the regularity of the solution map with respect to the random variables. Compared to other sampling based UQ methods, such as the Monte Carlo algorithm, SCMs will require far fewer realizations in order to provide accurate results. A prominent feature of these schemes is their {\it nonintrusive} nature, which makes it easy to combine them with any legacy PDE solver.

A recent approach for cutting down on the computational cost of UQ schemes is reduced order models (ROMs). These methods seek a linear space which approximates the solution manifold well. At the discrete level, this results in the need to solve orders of magnitude smaller linear systems than \textcolor{black}{those required} by full order models.

In order to devise an efficient UQ algorithm to address these challenges, we propose combining  {\it nonintrusive} SCMs with new {\it nonintrusively} stabilized ROM schemes. A major advantage of this approach is that it can be easily used in conjunction with legacy flow solvers. Additionally, the nonintrusive ROM stabilization will significantly decrease the computational cost of the flow solver, while still be able to tackle flows in the convection-dominated regime.

Our framework has several novel components: To our knowledge, SCM-ROM for the time-dependent NSE is new.  
SCM-ROM has been previously used only for the steady-state NSE~\cite{elman2013reduced}.
In fact, standard (i.e., without a ROM component) SCM investigations of parabolic PDEs are rather scarce:
SCM have been used for the heat equation~\cite{zhang2012error}, the convection-diffusion equation~\cite{GS07} (see also~\cite{torlo2017stabilized} for a reduced basis method approach), and the time-dependent NSE~\cite{tran2014convergence}. 

The second novel component of our framework is the nonintrusive filter-based stabilization of the SCM-ROM.
Specifically, at each timestep, given the flow variables at the current timestep $(u^{n},p^{n})$, we use the following {\it evolve-filter-relax (EFR)} algorithm to find the flow variables at the new timestep:
(1) evolve the flow variables one timestep using the given CFD code and call these new variables $(\tilde{u}^{n},\tilde{p}^{n})$;  
(2) use a spatial filter $F$ to filter the flow variables obtained in step (1): $(u^{n+1},p^{n+1}) = (F (\tilde{u}^{n}), F(\tilde{p}^{n}))$;
(3) relax the numerical approximation obtained in (2).
Using the spatial filter in step (2) decreases the numerical oscillations of the flow variables and, therefore, increases the numerical stability of the algorithm.
The role of the relaxation step (3) is to control the amount of numerical dissipation introduced in step (2).
The main advantage of the EFR algorithm over other stabilization algorithms is that it is {\it nonintrusive}, i.e., one can use the EFR algorithm {\it without modifying the CFD solver}.
Thus, the EFR algorithm can be easily used in legacy codes:
One can simply start with a legacy CFD code and, with minimal coding effort, can produce a code that can tackle moderate to high Reynolds number flows. 
Its nonintrusive character and simplicity have made the EFR algorithm extremely appealing for the numerical simulation of convection-dominated flows with standard methods, e.g., spectral methods~\cite{boyd1998two,boyd2001chebyshev,CHQZ88,hesthaven2008filtering,pasquetti2002comments}, spectral element methods~\cite{fischer2001filter,mullen1999filtering}, and finite element methods~\cite{ervin2012numerical,layton2014numerical}. 
However, we emphasize that the EFR algorithm has never been used in conjunction with ROMs; only an EF algorithm (without relaxation) has been used in~\cite{wells2017evolve}.
In this paper, we build on the initial efforts made in~\cite{wells2017evolve} and propose two new developments:
(i) we add a {\it relaxation} step to the EF-ROM introduced in~\cite{wells2017evolve}; and
(ii) we propose a novel {\it higher-order ROM differential filter}.

\medskip

The new stabilized SCM-ROM framework proposed above is an {\it efficient, practical solution for the UQ of realistic, convection-dominated fluid flows with legacy codes}.
The new SCM-ROM tackles the three challenges mentioned in the beginning of this section:
(i) It is computationally efficient; the computational cost of the ROM component is orders of magnitude lower than the cost of a brute-force DNS.
(ii) Its filter-based stabilization allows the numerical simulation of realistic, convection-dominated flows.
(iii) Both the SCM and ROM stabilization components are nonintrusive, which allows the new framework to be successfully used in conjunction with legacy codes.
We note that legacy codes are often used in the DNS of realistic flows and are consistently becoming more common in the ROM community (see, e.g., RBniCS~\cite{hesthaven2015certified} and pyMOR~\cite{milk2016pymor}).
Thus, nonintrusive discretizations are highly desirable for both DNS and ROMs of realistic flows.

\medskip

The rest of the paper is organized as follows:
In Section~\ref{sec:rom-uq}, we introduce the problem setting and present a brief overview of SCMs. 
In Section~\ref{sec:rom}, we briefly present ROMs for fluid flows.
In Section~\ref{sec:stabilized-scm-rom}, we propose stabilized SCM-ROMs for the numerical investigation of UQ in the time-dependent NSE in the convection-dominated regime.
In Section~\ref{sec:numerical-results}, we present numerical results for the new stabilized SCM-ROM for the two-dimensional (2D) flow past a circular cylinder with a random viscosity that yields a \textcolor{black}{random} Reynolds number with mean $Re=100$.
Finally, we draw conclusions in Section~\ref{sec:conclusions}.

\section{SCMs for Fluid Flows}
\label{sec:rom-uq}

\subsection{Problem Setting}

Let $D \subset \mathbb{R}^{s}$, $s= 2,3$, be an open regular domain with boundary $\partial D$ and let $[0,T]$ be a given time interval. Denote by $(\Omega, \mathcal{F}, P)$ a complete probability space, with $\Omega$ the set of outcomes, $\mathcal{F} \subset 2^{\Omega}$ the $\sigma$-algebra of events, and $P: F \rightarrow [0,1]$ a probability measure. The viscosity $\nu(x,\omega)$ is modeled as a random field, with $\omega \in \Omega$. Using no-slip boundary conditions, we consider the stochastic time-dependent Navier-Stokes equations:    
\begin{equation}\label{eq:NSE}
\left\{\begin{aligned}
u_{t} +u\cdot\nabla u -\nu(x,\omega) \Delta u +\nabla p  &
=f(x,t)&\quad\forall x\in D \times(0,T]\\
\nabla\cdot u  &  =0&\quad\forall x\in\ D \times(0,T]\\
u  &  =0&\quad\forall x\in\partial D \times(0,T]\\
u(x,0,\omega)  &  =u_{0}(x)&\quad\forall x\in D.
\end{aligned}\right.
\end{equation}
Here $u$ and $p$ denote the velocity and pressure of the flow, respectively\textcolor{black}{, and $f$ is the body force}.

{This problem can be put into a corresponding weak formulation. Defining the Hilbert spaces for velocity $X = H_{0}^{1}(D)$, pressure $Q = L_{0}^{2}(D)$, and stochastic space $W = L^{2}_{P}(\Omega)$, the weak formulation of~\eqref{eq:NSE} can be written as: find $u \in X \otimes W$ and $p \in Q \otimes W$ which, for almost all $t\in(0,T]$, satisfy 
	\begin{equation}\label{wfwf}
	\left\{\begin{aligned}
	\mathbb{E}[(u_{t},v)]+ \mathbb{E}[(u\cdot\nabla u,v)]+ \mathbb{E}[\nu(\nabla u,\nabla v)]- \mathbb{E}[(p
	,\nabla\cdot v)]  &  = \mathbb{E}[(f,v)]&\quad\forall v\in X \otimes W\\
	\mathbb{E}[(\nabla\cdot u,q)]  &  =0&\quad\forall q\in Q \otimes W\\
	u(x,0,\omega)  &  = u_{0}(x).&
	\end{aligned}\right.
	\end{equation}

	In a number of settings it is valid to assume that the randomness  can be approximated well by a finite number of random variables. Letting $y = (y_1,\ldots,y_d) \in \Gamma \subset \mathbb{R}^{d}$ be a finite $d$ dimensional vector distributed according to a joint probability density function $\rho(y)$ in some parameter space $\Gamma = \prod_{\ell=1}^{d} \Gamma_{\ell}$, the random field $\nu(x,\omega)$ can then be written in terms of the finite random variable as $\nu(x,y)$. Denoting by $Y = L^{2}_{\rho}(\Gamma)$ the space of square integrable functions on $\Gamma$ with respect to the weight $\rho(y)$, the weak formulation we will consider throughout the rest of this paper is: find $u \in X \otimes Y$ and $p \in Q \otimes Y$ which, for almost all $t\in(0,T]$, satisfy
	\begin{equation}\label{weak_formulation_final}
	\left\{\begin{aligned}
	\int_{\Gamma} \int_D (u_{t},v)\rho(y)dy + \int_{\Gamma} \int_D(u\cdot\nabla u,v)\rho(y)dy &+ \int_{\Gamma} \int_D \nu(x,y)(\nabla u,\nabla v)\rho(y)dy \  \\  -\int_{\Gamma} \int_D (p ,\nabla\cdot v)\rho(y)dy 
	&= \int_{\Gamma} \int_D (f,v)\rho(y)dy \quad\forall v\in X \otimes Y\\
	\int_{\Gamma} \int_D (\nabla\cdot u,q) \rho(y)dy  &  =0 \qquad \qquad \qquad \qquad \forall q\in Q \otimes Y\\
	u(x,0,\omega)   &=u_{0}(x).
	\end{aligned}\right.
	\end{equation}

	In order to ensure an efficient ROM scheme, we make an assumption of affine dependence on the random variable for the viscosity, i.e., 
	\begin{equation}
	\nu(x,y) = \nu_{0}(x) + \sum_{\ell=1}^{d} \nu_{\ell}(x)y_{\ell}.
	\end{equation}
	In the case that this assumption does not hold, empirical interpolation schemes \cite{CS10} have been devised to provide approximations that do fulfill this assumption.}

\subsection{Stochastic Collocation Methods (SCMs)}
\label{sec:SCM}

SCMs, which take advantage of the regularity of the solution map with respect to the random variables,  have recently been developed for the UQ of quantities of interest (QOI) of stochastic PDE solutions, $\psi(u)$. Some common choices for $\psi$ in the case of fluids include the drag, lift, and energy. The goal of SCMs is to obtain statistical information about this QOI, i.e., the expectation
\begin{equation} 
\label{QOI_expecation}
\mathbb{E}[\psi(u)] = \int_{\Gamma} \psi(u,y)\rho(y)dy.
\end{equation}

Two of the most popular SCMs are the  discrete least squares collocation method and the sparse grid algorithm \cite{CCMFT2013,lee2018stochastic,FTW2008}. In either of these schemes, for a given time $t_n$ and a set of sample points $y_{sc}$ = $\{y^{j}\}_{j=1}^{N_{sc}} \subset \Gamma$, the PDE is solved using some spatial approximation, giving the solution denoted by $u(x,t_n,y^{j})$. Then, given a basis $\{\phi_{\ell}\}$ for the space $L^{2}_{\rho}(\Gamma)$, a discrete approximation of the form 
\begin{equation}
\label{SC_approx}
u^{sc} (x,t_n,y) := \sum_{\ell=1}^{d} c_{\ell}(x)\phi_{\ell}(y)
\end{equation}
is constructed with the coefficients $c_{\ell}(x)$.

There is a wide range of options for both the approximation points $y_{sc}$ and polynomial basis $\{\phi_{\ell}\}$. A common choice for the basis functions will be polynomials which are orthogonal with respect to the PDF $\rho$. The sample points will depend upon the SCM used. For discrete least squares algorithms, common approaches are to take random samples with respect to $\rho$. When a sparse grid approach is taken, \eqref{SC_approx} becomes an interpolant. Some commonly used interpolations points include nested rules, such as Leja and Clenshaw Curtis points. An advantage of these interpolation schemes is that they have associated quadrature weights $\{w^{j} \}_{j=1}^{N_{sc}}$. In this case, an approximation to \eqref{QOI_expecation} will be given as 
\begin{equation}
\label{SG_quad}
\mathbb{E}[\psi(u)] = \int_{\Gamma} \psi(u,y)\rho(y)dy \approx \sum_{j = 1}^{N_{sc}} w^{j} \psi(u,y^{j}).
\end{equation}

\subsection{SCMs for Fluid Flows}
\label{sec:scm-fluid-flows}

To our knowledge, there are relatively few investigations of SCMs for the NSE:
In \cite{elman2013reduced}, SCMs were combined with a reduced basis method and applied to the steady NSE. 
In \cite{Sankaran2011ASC}, SCMs were applied to the time-dependent NSE and fast convergence was demonstrated for short term time evolutions. 
In \cite{WK06}, the long term behavior of SCMs for time-dependent NSE was studied using a polynomial chaos method for random frequency stochastic processes. 
Finally, the numerical investigation in~\cite{WK06} was generalized in \cite{tran2014convergence} to a wider range of stochastic processes; theoretical foundations for the numerical results were also proven in \cite{tran2014convergence}.

\section{ROMs for Fluid Flows}
\label{sec:rom}

In this section, the standard Galerkin ROM for fluid flows is presented succinctly.

\subsection{ROM Basis}
\label{sec:pod}

As a ROM basis, we use the proper orthogonal decomposition (POD)~\cite{HLB96,noack2011reduced,Sir87abc}, which we briefly describe next.  
We emphasize that the new framework could be extended to other ROM bases, e.g., the dynamic mode decomposition (DMD)~\cite{kutz2016dynamic,rowley2009spectral,schmid2010dynamic}.
The POD starts with the snapshots $\{u^1_h,\ldots, u^{N}_h\}$, which are numerical approximations of~\eqref{eq:NSE} at $N$ different time instances.
The finite element (FE) solutions of~\eqref{eq:NSE} are considered as snapshots in this section.
We emphasize, however, that other numerical methods can be used instead.
The POD seeks a low-dimensional basis that approximates the snapshots
optimally with respect to a certain norm. In this paper, the commonly
used $L^2$-norm will be chosen. The solution of the minimization
problem is equivalent to the solution of the eigenvalue problem
$
\mathbb{A}\mathbb{A}^\top M_{h} \varphi_j = \lambda_j \varphi_j,
\  j=1,\ldots,\mathcal{N},
$
where $\varphi_j$ and $\lambda_j$ denote the vector of the FE coefficients
of the POD basis functions and the POD eigenvalues, respectively, $\mathbb{A}$
denotes the snapshot matrix, whose columns correspond to the FE
coefficients of the snapshots, $M_{h}$ denotes the FE mass matrix, and $\mathcal{N}$
is the dimension of the FE space $X^h$.
The eigenvalues are real and non-negative, so they can be ordered as
follows:
$
\lambda_1 \ge \lambda_2 \ge \ldots \ge \lambda_R > \lambda_{R + 1}
= \ldots = \lambda_{\mathcal{N}} = 0.
$
The POD basis consists of the normalized functions $\{
\varphi_{j}\}_{j=1}^{r}$, which correspond to the first $r\le \mathcal{N}$ largest
eigenvalues. Thus, the POD space is defined as $X^r := \text{span} \{
\varphi_1, \ldots, \varphi_r \}$.

\subsection{Galerkin ROM (G-ROM) for Fluid Flows}
\label{sec:g-rom}

To develop the standard Galerkin ROM for fluid flows, we start by considering the POD approximation of the velocity
\begin{equation}
{u}_r(x,t)
\equiv 
\sum_{j=1}^r a_j(t) \varphi_j(x),
\label{eqn:g-rom-1}
\end{equation}
where $\{a_{j}(t)\}_{j=1}^{r}$ are the sought time-varying coefficients that represent the POD-Galerkin trajectories.
Since the POD basis functions are simply linear combinations of the snapshots, we assume that each POD basis function is weakly divergence-free and the compressibility constraint is automatically satisfied (see, however,~\cite{hesthaven2015certified,quarteroni2015reduced} for alternatives).
By using the POD basis in a Galerkin approximation of the NSE, the standard {\it Galerkin ROM (G-ROM)} is obtained: $\forall \, i = 1, \ldots, r,$
\begin{eqnarray}
\left( \frac{\partial u_{r}}{\partial t} , \varphi_i \right)
+ \nu \, \biggl( \nabla u_{r} , \nabla \varphi_i \biggr)
+ \biggl( (u_{r} \cdot \nabla) u_{r} , \varphi_i \biggr)
= 0 \, .
\label{eqn:g-rom}
\end{eqnarray}
The G-ROM~\eqref{eqn:g-rom} yields the following autonomous dynamical system for the vector of time coefficients, $a(t)$:
\begin{equation}
\dot{a} = A a   + a^\top B a ,
\label{eqn:g-rom-3}
\end{equation}
\noindent
where $A$ and $B$ correspond to the linear and quadratic terms in the numerical discretization of the NSE~\eqref{eq:NSE}, respectively.
The finite dimensional system \eqref{eqn:g-rom-3} can be written componentwise  as follows:
For all $i = 1, \ldots, r$,
\begin{eqnarray}
\dot{a}_i(t)
=  
\sum_{m=1}^{r} A_{im}a_m(t) 
+ \sum_{m=1}^r \sum_{n=1}^r B_{imn}a_n(t)a_m(t),
\label{eqn:g-rom-5}
\end{eqnarray}
where
$
A_{im}
= 
- \nu \left( \nabla \varphi_m , \nabla \varphi_i \right)
$
and
$
B_{imn}
= - \bigl( \varphi_m \cdot \nabla \varphi_n , \varphi_i \bigr) \, .
$

\subsection{Stabilized ROMs for Convection-Dominated Flows}

The numerical discretization of convection-dominated flows with standard Galerkin methods (e.g., FEs) generally displays spurious numerical oscillations~\cite{john2016finite}.
Thus, stabilized numerical discretizations are used instead~\cite{john2016finite}.

Since G-ROM~\eqref{eqn:g-rom} is a Galerkin method, it also yields numerical oscillations when used in the numerical simulation of convection-dominated flows~\cite{HLB96}.
Thus, many stabilized ROMs have been proposed over the years; see, e.g.,~\cite{balajewicz2016minimal,ballarin2015supremizer,benosman2017learning,bergmann2009enablers,fick2017reduced,gunzburger2018leray,rebollo2017certified,stabile2018finite,wells2017evolve} for recent work on this topic.
Among these stabilization methods, we mention the filter based stabilization proposed in~\cite{wells2017evolve} (see also~\cite{gunzburger2018leray,iliescu2018regularized,sabetghadam2012alpha, xie2018evolve} for related work), in which spatial filtering is used to smooth flow variables and alleviate the G-ROM's numerical instability. 
In the next section, we propose a new filter-based stabilization for convection-dominated flows.

\section{Evolve-Filter-Relax SCM-ROM}
\label{sec:stabilized-scm-rom}

In this section, we propose a new stabilized SCM-ROM for the numerical investigation of UQ in convection-dominated flows.
In Section~\ref{sec:rom-filters}, we first review the current spatial ROM filters and then propose a novel ROM spatial filter: The higher-order ROM differential filter.
In Section~\ref{sec:efr-rom}, we use the higher-order ROM differential filter to develop a new stabilized ROM, the evolve-filter-relax ROM.
Finally, in Section~\ref{sec:stabilized-scm-rom}.3, we combine the new evolve-filter-relax ROM with the SCM to obtain a stabilized SCM-ROM for convection-dominated flows.

\subsection{ROM Spatial Filters}
\label{sec:rom-filters}

To build the evolve-filter-relax ROM in Section~\ref{sec:efr-rom}, we use explicit ROM spatial filtering to smooth the flow variables and increase the numerical stability of the models.
Thus, we need to address the following problem:
Given the current ROM approximation at timestep $n$
\begin{eqnarray}
u_{r}^{n}
= \sum_{j=1}^{n} a_{j}^{r} \varphi_{j} \, ,
\label{eqn:rom-df-0}
\end{eqnarray}
define the filtered ROM approximation at timestep $n$, $\overline{u_{r}^{n}}$.

In this section, we present the ROM differential filter~\cite{wells2017evolve,xie2017approximate} (Section~\ref{sec:df}) and propose a {\it novel higher-order ROM differential filter} (Section~\ref{sec:hodf-rom}); see, e.g.,~\cite{xie2018data} for alternative ROM spatial filters.

\subsubsection{Differential Filter (DF)}
\label{sec:df}

For a given $u_{r}^{n} \in X^r$ and a given radius $\delta$, the {\it differential filter (DF)} computes $\overline{u_{r}^{n}} \in X^r$ as follows:
\begin{eqnarray}
\biggl(
\left(
I - \delta^2 \Delta
\right) \overline{u_{r}^{n}} , \varphi_j
\biggr)
= (u_{r}^{n}, \varphi_j),
\quad \forall \, j=1, \ldots r \, .
\label{eqn:df}
\end{eqnarray}
We emphasize that the DF uses an {\it explicit length scale} $\delta$ (i.e., the radius of the filter).
Differential filters have been used in the simulation of convection-dominated flows with standard numerical methods~\cite{germano1986differential-b,germano1986differential}.
In reduced order modeling, the DF was used in~\cite{gunzburger2018leray,iliescu2018regularized,sabetghadam2012alpha,wells2017evolve,xie2017approximate}.

The DF~\eqref{eqn:df} has several appealing properties:
(i) It acts as a spatial filter, since it eliminates the small scales (i.e., high frequencies) from the input. 
Indeed, the DF~\eqref{eqn:df} uses an elliptic operator to smooth the input variable. 
(ii) The DF has a low computational overhead, since it amounts to either solving a linear system with a very small $r \times r$ matrix that is precomputed or precomputing and storing the filtered ROM modes.
(iii) The DF preserves incompressibility in the NSE, since it is a linear operator.

There are two types of DF discretization that are currently in use: A ROM version and an FE version~\cite{wells2015phd}.
In our numerical investigation, we use the ROM discretization of the DF; we emphasize, however, that the new stabilized ROM that we propose can \textcolor{black}{also} be used with the FE discretization of the DF.
The ROM discretization of the DF is defined as follows:
Let $\delta$ be the radius of the ROM-DF.
For all $u_r \in X^{r}$, find $\orr{u_r} \in X^r$ such that
\begin{eqnarray}
- \delta^2 \, \Delta \orr{u_{r}^{n}}
+ \orr{u_{r}^{n}}
&=& u_{r}^{n}
\qquad \text{in } \Omega
\label{eqn:rom-df-1} \\
\orr{u_{r}^{n}} 
&=& 0  
\qquad \text{on } \partial \Omega \, ,
\label{eqn:rom-df-2}
\end{eqnarray}
which yields the following linear system:
\begin{eqnarray}
\left( M_r + \delta^2 \, S_r \right) \, \orr{a_r^{n}}
= a_r^{n} \, ,
\label{eqn:pod-df-linear-system}
\end{eqnarray}
where 
$M_r \in \R^{r \times r}$ is the \textcolor{black}{ROM} mass matrix with entries $(M_r)_{i j} = ( \varphi_j , \varphi_i)$,
$S_r \in \R^{r \times r}$ is the \textcolor{black}{ROM} stiffness matrix with entries $(S_r)_{i j} = (\nabla \varphi_j , \nabla \varphi_i)$,
$\orr{a_r^{n}} \in \R^{r}$ is the vector of \textcolor{black}{ROM} coefficients of the output filtered variable $\orr{u_{r}^{n}}$, and
$a_r^{n} \in \R^{r}$ is the vector of \textcolor{black}{ROM} coefficients of the input variable $u_{r}^{n}$.

We emphasize that the computational cost of the ROM-DF is much lower than the computational cost of the FE-DF, since the former yields a very small $r \times r$ linear system, whereas the latter yields a much larger linear system.

\subsubsection{Higher-Order Differential Filter (HODF)}
\label{sec:hodf-rom}

In this section, we propose a novel {\it higher-order ROM differential filter (HODF)}.

The need for higher-order filters has been emphasized in~\cite{hesthaven2008filtering} for spectral methods,  in~\cite{mullen1999filtering} for SEM, and in~\cite{san2015posteriori} for finite difference methods.
For example, the left panel of Figure~1 in~\cite{mullen1999filtering} shows that the transfer function of the DF does not have a very sharp decay, whereas the higher-order filters in the right panel of Figure~1 do have a sharp decay that approaches the decay of the sharp cut-off filter (see also Figures. 4-6 in~\cite{san2015posteriori}).

In this section, we extend to ROM arena the high-order SEM filters proposed in~\cite{mullen1999filtering}.

\paragraph{HODF -- Version 1:}
The easiest way to enforce a sharp decay of the transfer function of the DF is to apply repeatedly the DF filter (see left panel of Figure~1 in~\cite{mullen1999filtering}).
This amounts to replacing~\eqref{eqn:pod-df-linear-system} with the following:
\begin{eqnarray}
\left( M_r + \delta^2 \, S_r \right)^{m} \, \orr{a_r^{n}}
= a_r^{n} \, ,
\label{eqn:hodf-1}
\end{eqnarray}
where $m$ is the order of the HODF.
The major criticism in~\cite{mullen1999filtering} was that the transfer function of the HODF~\eqref{eqn:hodf-1} does not have a sharp decay (see left panel of Figure~1 in~\cite{mullen1999filtering}).

\paragraph{HODF -- Version 2:}
To enforce an even sharper decay of the transfer function of the HODF~\eqref{eqn:hodf-1}, we propose the following:
\begin{eqnarray}
\left( M_r + \delta^2 \, S_r^{m} \right) \, \orr{a_r^{n}}
= a_r^{n} \, ,
\label{eqn:hodf-2}
\end{eqnarray}
where $m$ is the order of the HODF.
This effectively amounts to replacing the Laplacian $\Delta$ in the standard DF~\eqref{eqn:df} with an $m$-th order Laplacian $\Delta^{m}$.

The transfer function of the HODF~\eqref{eqn:hodf-2} is significantly sharper than the transfer function of the HODF~\eqref{eqn:hodf-1} (compare left and right panels of Figure~1 in~\cite{mullen1999filtering} and also Figures 3-6 in~\cite{san2015posteriori}).
However, the authors in~\cite{mullen1999filtering} emphasize that the condition number of $S_r^{m}$ is extremely high.  
Thus, they propose to use a modified CG algorithm.

\begin{remark}[Previous Relevant Work]
	The first HODF was proposed in~\cite{fischer2001filter} for spectral element methods.
	For finite difference methods, HODFs were investigated in~\cite{san2015posteriori}.
	To our knowledge, the HODFs in~\eqref{eqn:hodf-1} and \eqref{eqn:hodf-2} are novel in a ROM setting.
\end{remark}

\begin{remark}[Computational Cost]
	We emphasize that the computational cost of the HODF (and DF, for that matter) is negligible (at least for their ROM versions), since this filter involves only low-dimensional ($\mathcal{O}(10)$) matrices.
	This is in stark contrast with the SEM setting, where HODF's computational cost is prohibitively high.
\end{remark}

\subsection{Evolve-Filter-Relax ROM (EFR-ROM)}
\label{sec:efr-rom}

In this section, we put forth an {\it evolve-filter-relax ROM (EFR-ROM)} as a {\it nonintrusive} stabilization method for ROMs of convection-dominated flows.
We use the DF and HODF filters proposed in Section~\ref{sec:rom-filters} as ROM spatial filters for the EFR-ROM.
Since the computational cost of the DF and HODF is negligible, the computational cost of the EFR-ROM will be very low, similar to the computational cost of the G-ROM.
The EFR-ROM's nonintrusive character and simplicity allow one to simply start with a legacy CFD/ROM code and, with minimal coding effort, produce a ROM code that can tackle moderate to high Reynolds number flows. 

The EFR algorithm has been used in conjunction with classical numerical methods, e.g., spectral methods~\cite{boyd1998two,boyd2001chebyshev,CHQZ88,pasquetti2002comments}, spectral element methods~\cite{fischer2001filter,mullen1999filtering}, finite difference methods~\cite{mathew2003explicit}, 
and finite element methods~\cite{ervin2012numerical,layton2014numerical}.
To our knowledge, the EFR model has not been used in a ROM setting.
Indeed, the only current evolve-filter ROM (EF-ROM) is that proposed in~\cite{wells2017evolve}; this EF-ROM, however, did not include a relaxation step and, as such, was overly-dissipative.

\medskip

The algorithm for the new EFR-ROM reads: $\forall \, n = 1, \ldots, N$ and $\forall \, i = 1, \ldots, r,$

\begin{eqnarray}
&&	\text{\bf Evolve:}
\nonumber \\[0.2cm]
&& \quad 
\left(
\frac{3 w_{r}^{n+1} -4 u_{r}^{n}+ u_{r}^{n-1}}{2\Delta t} , \varphi_{i}
\right)
+ \nu \, \biggl(
\nabla w_{r}^{n+1} ,
\nabla \varphi_{i}
\biggr)
+ \biggl(
(2u_{r}^{n}-u_{r}^{n-1}) \cdot \nabla w_{r}^{n+1} ,
\varphi_{i}
\biggr)
\nonumber \\[0.2cm]
&& \hspace*{10.0cm} 
= \biggl(
f^{n+1} ,
\varphi_{i}
\biggr) \, ,
\label{eqn:ef-rom-1} \\[0.3cm]
&&	\text{\bf Filter:} \qquad 
\text{Compute} \quad   
\overline{w_{r}^{n+1}} \, ,
\label{eqn:ef-rom-2} \\[0.3cm]
&&	\text{\bf Relax:} \qquad 
u_{r}^{n+1}
= (1 - \chi) \, w_{r}^{n+1}
+ \chi \, \overline{w_{r}^{n+1}} \, ,
\label{eqn:ef-rom-3}
\end{eqnarray}
\vspace*{-0.3cm}

\noindent where $\Delta t$ is the timestep and $\chi \in [0,1]$ is the relaxation parameter.
The ``evolve'' step in the EFR-ROM (i.e., equation~\eqref{eqn:ef-rom-1}) is just one step of the time discretization of the standard G-ROM~\eqref{eqn:g-rom}.
The ``filter'' step in the EFR-ROM (i.e., equation~\eqref{eqn:ef-rom-2}) consists of filtering of the intermediate solution obtained in the ``evolve'' step with the DF and HODF proposed in Section~\ref{sec:rom-filters}.
Finally, the ``relax" step in the EFR-ROM (i.e., equation~\eqref{eqn:ef-rom-3}) averages the unfiltered and filtered flow variables $w_{r}^{n+1}$ and $\overline{w_{r}^{n+1}}$, respectively.

We note that a BDF2 time discretization was used in the EFR-ROM~\eqref{eqn:ef-rom-1}--\eqref{eqn:ef-rom-3}, but other time discretizations are possible~\cite{ervin2012numerical}.

As pointed out in~\cite{ervin2012numerical}, if $\chi=1$ (i.e., there is no relaxation) and an explicit time-discretization is used, the EFR-ROM used with the DF 
reduces to the standard G-ROM plus numerical diffusion of magnitude $\frac{\delta^2}{\Delta t}$.
(We note, however, that if an implicit time discretization is used or if one does not use the DF, the relaxation effect is not that clear~\cite{ervin2012numerical}.
In fact, the time-relaxation seems to be quite different from the DF; Leo Rebholz, personal communication.)
Thus, to diminish the magnitude of the numerical diffusion, an extra relaxation step is generally used, i.e., $\chi<1$~\cite{ervin2012numerical,fischer2001filter}.
The relaxation step decreases the amount of numerical diffusion and increases the accuracy; see, e.g., the numerical results in~\cite{bertagna2016deconvolution,fischer2001filter} and the theoretical results in~\cite{ervin2012numerical}.
The scaling $\chi \sim \Delta t$ is commonly used in EFR models~\cite{ervin2012numerical}.
In~\cite{bertagna2016deconvolution}, however, the authors advocate higher values for $\chi$.

\begin{remark}
	We emphasize that the EFR-ROM is fundamentally different from the {\it preprocessing} utilized to filter out the noise in the snapshot data~\cite{aradag2011filtered} (which is a G-ROM that uses the filtered basis functions, $\overline{\varphi_{i}}$).
	The reason is that the EFR-ROM, spatial filtering is performed at each timestep, which amounts to filtering the ROM basis functions repeatedly. 
	Since $\overline{\varphi} \neq \overline{\overline{\varphi}}$, we conclude that EFR-ROM is different from the preprocessing approach.
\end{remark}

\subsection{EFR-SCM-ROM}
\label{sec:efr-scm-rom}
In this section, we provide a detailed description of the combination of the EFR-ROM scheme given in Section \ref{sec:efr-rom} with the SCMs described in Section \ref{sec:SCM}. The EFR-SCM-ROM algorithm can be split into an offline and online phase. The overall goal of the EFR-SCM-ROM is to evaluate the approximation $u_h(x,t,y)$ at a sufficiently large number of times $t_n$ and sample points $N_{sc}$ generated by a chosen SCM scheme so that the approximation \eqref{SC_approx} is sufficiently accurate. However, due to a lack of computational resources, we will only be able to obtain full order solution trajectories for a smaller number of sample points, $N_{train}$. In the offline phase, we generate discrete solutions (known as snapshots) for these $N_{train}$ sample points. Using POD, we build a reduced basis from this snapshot set, which should accurately approximate the data present in the snapshot set. In the online phase, the POD basis will now be used to generate approximations for other parameter values determined by the targeted $N_{sc}$ collocation points.

In the rest of this subsection, we provide details on each phase, as well as a full algorithm for the EFR-SCM-ROM. For other examples of SCM and ROM combinations, see \cite{elman2013reduced,chen2014comparison}.

\subsubsection{Offline Phase} 
In the offline phase, we are interested in constructing a reduced basis in both time and stochastic space. Let $y_{sc}$ = $\{y^{j}\}_{j=1}^{N_{sc}}$ denote the $N_{sc}$ total stochastic collocation points we will compute with in the online phase and let $y_{train}$ = $\{y^{j}\}_{j=1}^{N_{train}}$ denote the training set which we will use in the offline phase (where $N_{train} << N_{sc}$). Inputting $y_{train}$ into the parameterized viscosity will result in a total of $\nu_{train}$ = $\{\nu_{k}\}_{k=1}^{N_{train}}$ viscosities. For a positive integer $N'$, we then let $0=t_0<t_1< \cdots < t_m < \cdots < t_{N'} = T_{train}$ denote a uniform partition of the time interval $[0,T_{train}]$, with $T_{train} \leq T$. Denoting by $\vec{u}_S^{k,m}$ a finite dimensional solution vector to the NSE (i.e., FE solution) at time $t_m$ for viscosity $\nu_k$, we construct the snapshot matrix 
$$
\mathbb{A} = \big(\vec{u}_S^{1,0},\vec{u}_S^{1,1}, \ldots , \vec{u}_S^{1,N_T}, \vec{u}_S^{2,0},\vec{u}_S^{2,1},  \ldots , \vec{u}_S^{2,N_T}, \ldots , \vec{u}_S^{N_{train},0},\vec{u}_S^{N_{train},1}, \ldots , \vec{u}_S^{N_{train},N'}\big).
$$ 
Solving the constrained minimization problem outlined in Section \ref{sec:pod}, we then construct the reduced basis $\{\varphi_{j}\}_{j=1}^{r}$.

\subsubsection{Online Phase}
In the online phase, we compute the approximation \eqref{SC_approx} using the targeted sample points $\{y^{j}\}_{j=1}^{N_{sc}}$. Inputting $y_{train}$ into the parameterized viscosity will result in a total of $\nu_{sc}$ = $\{\nu_{k}\}_{k=1}^{N_{sc}}$ viscosities. We now consider the time interval $0=t_0<t_1< \cdots < t_m < \cdots < t_{N} = T$. Fixing the parameters $\delta$ and $\chi$ for the filter and relaxation step respectively, we compute solutions to the NSE at each timestep $u_{r,k}^{n+1}$ using the EFR scheme for each individual viscosity. This reduces to solving at each timestep  
$\forall \, k = 1, \ldots, N_{sc}$ and $\forall \, i = 1, \ldots, r,$

\begin{eqnarray}
&&	\text{\bf Evolve:}
\nonumber \\[0.2cm]
&& \quad 
\left(
\frac{3 w_{r,k}^{n+1} -4 u_{r,k}^{n}+ u_{r,k}^{n-1}}{2\Delta t} , \varphi_{i}
\right)
+ \nu_{k} \, \biggl(
\nabla w_{r,k}^{n+1} ,
\nabla \varphi_{i}
\biggr)
+ \biggl(
(2u_{r,k}^{n}-u_{r,k}^{n-1}) \cdot \nabla w_{r,k}^{n+1} ,
\varphi_{i}
\biggr)
\nonumber \\[0.2cm]
&& \hspace*{10.0cm} 
= \biggl(
f^{n+1} ,
\varphi_{i}
\biggr) \, ,
\label{eqn:ef-rom-1_SC} \\[0.3cm]
&&	\text{\bf Filter:} \qquad 
\text{Compute} \quad   
\overline{w_{r,k}^{n+1}} \, ,
\label{eqn:ef-rom-2_SC} \\[0.3cm]
&&	\text{\bf Relax:} \qquad 
u_{r,k}^{n+1}
= (1 - \chi) \, w_{r,k}^{n+1}
+ \chi \, \overline{w_{r,k}^{n+1}} \, ,
\label{eqn:ef-rom-3_SC}
\end{eqnarray}

For each solution, we then compute the corresponding quantity of interest $\psi(u_{r,k}^{n+1})$. Applying the quadrature formula \eqref{SG_quad} gives an approximation to $\mathbb{E}[\psi(u^{n+1})]$. A full outline of the EFR-SCM-ROM scheme is given in Algorithm \ref{ERF-SCM-ROM-algo}.

\begin{algorithm}
	\caption{EFR-SCM-ROM}\label{ERF-SCM-ROM-algo}
	\begin{algorithmic}[1]
		\Procedure OFFLINE construction 
		\State \textbf{Initialization:} Mesh, FE functions, $T_{train}$, $N_{train}$, $u^{0}_h$, snapshot matrix $\mathbb{A}$, $\{y^{j}\}_{j=1}^{N_{train}}$, $ \{\nu_{k}\}_{k=1}^{N_{train}}$
		\For{$n = 0 \ldots N'$}
		\For{$k = 1 \ldots N_{train}$}
		\State Compute $u^{n+1}_{h,k}$
		\State Store: $\mathbb{A} = [\mathbb{A}, u^{n+1}_{h,k}]$
		\EndFor
		\EndFor
		\State Solve for the reduced basis $\{\varphi_j\}_{j=1}^{r}$ : $\mathbb{A}\mathbb{A}^\top M_{h} \varphi_j = \lambda_j \varphi_j$
		\  $j=1,\ldots,r$
		\EndProcedure
		\Procedure ONLINE Construction
		\State \textbf{Initialization:} $T$, $N_{sc}$, $u_{r}^{0}$,
		$\{y^{j}\}_{j=1}^{N_{sc}}$,$\{w^{j}\}_{j=1}^{N_{sc}}$, $ \{\nu_{k}\}_{k=1}^{N_{sc}}$, $\delta$, $\chi$
		\State \textbf{Pre-compute:} $M_{r},S_r,A_r,B_r$ 
		\For{$n = 0 \ldots N$}
		\For{$k = 1 \ldots N_{sc}$}
		\State Evolve: $w_{r,k}^{n+1}$ according to \eqref{eqn:ef-rom-1_SC}
		\State Filter: $\overline{w_{r}^{n+1}}$ according to \eqref{eqn:ef-rom-2_SC}
		\State Relax: $u_{r}^{n+1}$ according to \eqref{eqn:ef-rom-3_SC}
		\State Calculate $\psi(u_{r,k}^{n+1})$
		\EndFor
		\State Estimate $\mathbb{E}[\psi(u^{n+1})] \approx \sum_{k = 1}^{N_{sc}} w^{k} \psi(u_{r,k}^{n+1})$
		\EndFor
		\EndProcedure
	\end{algorithmic}
\end{algorithm}

\section{Numerical Results}
\label{sec:numerical-results}

In this section, we perform a numerical investigation of the new SCM-EFR-ROM equipped with the DF~\eqref{eqn:df} and the HODF~\eqref{eqn:hodf-1} with $m=2, 3$, and $4$.
Specifically, we investigate whether the new EFR-ROM is more accurate than the standard G-ROM~\eqref{eqn:g-rom} in the SCM setting.
As a benchmark for our numerical investigation, we use the FE approximation. 
Thus, we compare four models, which we denote as follows: 
(i) DF-ROM, which is SCM-EFR-ROM with the DF~\eqref{eqn:df};
(ii) HODF-ROM-$m$, which is SCM-EFR-ROM with the HODF and $m=2, 3$, and $4$;
(iii) G-ROM, which is SCM-ROM~\eqref{eqn:g-rom}; and
(iv) DNS, which, with an abuse of notation, is the FE approximation used to generate the snapshots.
For the HODF-ROM-$m$ model, we investigate both HODF versions proposed in Section~\ref{sec:hodf-rom}:  version 1 (i.e., equation~\eqref{eqn:hodf-1}) and version 2 (i.e., equation~\eqref{eqn:hodf-2}).
We compare the four models (i)--(iv) in the numerical investigation of a 2D flow past a circular cylinder with a random viscosity that yields a \textcolor{black}{random} Reynolds number with mean $Re=100$.
We consider two cases for the random viscosity: a 1D constant viscosity case and a 5D variable viscosity case. 
To compare the four models, we use the accuracy of the time evolution of the energy coefficients; we note, however, that the time evolution of the lift and drag coefficients display a similar behavior.

The SCM used throughout this section will be a Clenshaw Curtis sparse grid generated via the software package TASMANIAN \cite{stoyanov2015tasmanian,doecode_6305}.

The rest of this section is organized as follows:
In Section~\ref{sec:test-problem-setup}, we describe the test problem setup.
In Section~\ref{sec:numerical-results-dns}, we present details of the DNS discretization and in Section~\ref{sec:numerical-results-roms} we outline the ROM construction.
In Section~\ref{sec:numerical-results-1d}, we present results for the numerical investigation of the four models for a constant 1D random viscosity case.
Finally, in Section~\ref{sec:numerical-results-5d}, we perform a similar investigation for a variable 5D random viscosity case.

\subsection{Test Problem Setup}
\label{sec:test-problem-setup}

\begin{figure}[h!]
	\begin{center}
		\includegraphics[width=0.7\textwidth,height=0.24\textwidth, trim=0 0 0 0, clip]{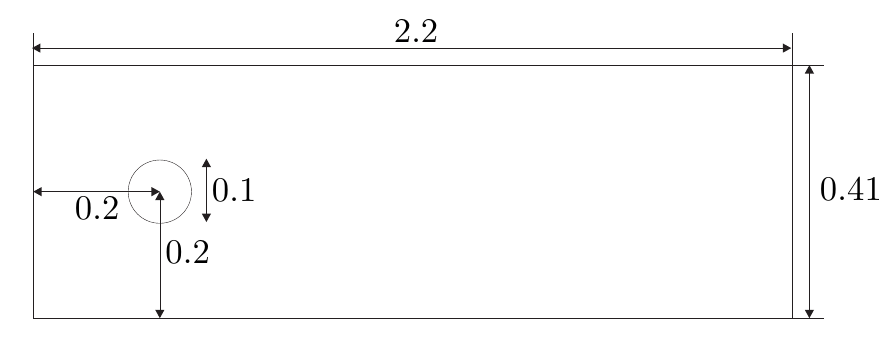}
	\end{center}
	\caption{\label{cyldomain} Channel flow around a cylinder domain.}
\end{figure}

The domain is a $2.2\times 0.41$ rectangular channel with a radius $=0.05$ cylinder, centered at $(0.2,0.2)$, see Figure~\ref{cyldomain}.  
No slip boundary conditions are prescribed for the walls and on the cylinder, and the inflow and outflow profiles are given by~\cite{mohebujjaman2018physically,xie2018data} $u_{1}(0,y,t)=u_{1}(2.2,y,t)=\frac{6}{0.41^{2}}y(0.41-y) \, , u_{2}(0,y,t)=u_{2}(2.2,y,t)=0$.
The kinematic viscosity $\nu(x,\omega)$ is random, there is no forcing, and the flow starts from rest. 
To compute the lift and drag coefficients, we use the following formulas~\cite{john2004reference,mohebujjaman2017energy}:
$$c_d(t)=\frac{2}{\tilde{\rho} L U_{max}^2}\int_{\partial S}\left(\tilde{\rho}\nu(x,y_i)\frac{\partial u_{ts}(t)}{\partial n}n_y-p(t)n_x\right)dS,$$

$$c_l(t)=\frac{2}{\tilde{\rho} L U_{max}^2}\int_{\partial S}\left(\tilde{\rho}\nu(x,y_i)\frac{\partial u_{ts}(t)}{\partial n}n_x-p(t)n_y\right)dS,$$
where $u_{ts}$ is the tangential velocity, $U_{max}=1$ is the maximum inlet velocity, $L=0.1$ is the cylinder diameter, $\tilde{\rho}=1$ is the density, $\partial S$ is the surface of the cylinder, and $<n_x,n_y>$ is the outward unit normal to the domain.

\subsection{DNS}
\label{sec:numerical-results-dns}

In this section, we present details of the algorithm \textcolor{black}{outlined} in Section~\ref{sec:efr-scm-rom} that is used to construct the benchmark DNS results.
To generate the DNS data, we run the time-dependent NSE~\eqref{eq:NSE} up to the simulation time $t=10$ for each of the collocation points $y_i\in\Omega$ and use the SCM to compute the quantity of interest $\psi$, i.e., the energy, lift, and drag coefficients. 
For the spatial discretization of the NSE, we use the $(X_h,Q_h)=(P_2,P_1^{disc})$ Scott-Vogelius (SV) element, which satisfies the mass conservation point-wise, on a \textcolor{black}{barycenter refined regular} triangular mesh providing $35020$ velocity and $25974$ pressure degrees of freedom. 
For the temporal discretization, we use a linearized BDF2 temporal discretization with a timestep size $\Delta t=0.002$. 
On the first timestep, we use a backward Euler temporal discretization to generate the \textcolor{black}{second} initial approximations required by the BDF2 scheme. 
The time discretization reads:
For $n=1,2,\cdots,$ find $(u_h^{n+1}, p_h^{n+1})\in (X_h, Q_h)$ satisfying for every $(v_h,q_h)\in(X_h,Q_h)$,
\begin{equation*}
\begin{aligned}
\left(\frac{3u_h^{n+1}-4u_h^n+u_h^{n-1}}{2\Delta t}, v_h\right)+((2u_h^{n}-u_h^{n-1})\cdot\nabla u_h^{n+1}, v_h)-(p_h^{n+1},\nabla\cdot v_h)\\+\nu(x,y_i)(\nabla u_h^{n+1},\nabla v_h) =0,\\
(\nabla\cdot u_h,q_h) = 0.
\end{aligned}
\end{equation*}

From the energy, drag, and lift plots in Figures~\ref{cons_vis_dns_energy}-\ref{cons_vis_dns_lift}, we observe that after an initial spin-up, the flow asymptotically reaches a statistically steady state by about $t=5$. 

\subsection{ROMs}
\label{sec:numerical-results-roms}

In this section, we present details of the EFR-SCM-ROM algorithm presented in Section~\ref{sec:efr-scm-rom}, which is used to construct the DF-ROM and HODF-ROM-$m$.
In Section~\ref{sec:numerical-results-roms-offline}, we present details of the offline phase and in Section~\ref{sec:numerical-results-roms-online}, we present details of the online phase.
In Section~\ref{sec:numerical-results-roms-computational-efficiency}, we comment on the computational efficiency of the DF-ROM and HODF-ROM-$m$.

\subsubsection{Offline Phase}
\label{sec:numerical-results-roms-offline}

In the offline phase of the EFR-SCM-ROM algorithm, we construct the DF-ROM and HODF-ROM-$m$ (see Section~\ref{sec:efr-scm-rom}). 

To this end, we generate a number of $N_{train}$ collocation points $\{y_i\}_{i=1}^{N_{train}}$, where $N_{train}<<N_{sc}$. 
We run our simulations for all these $N_{train}$ collocation points up to $t=10$. 
We collect a total of $2500$ snapshots at each timestep from $t=5$ to $t=10$ for each of the collocation points. 
We put together the resulting $2500\times N_{train}$ snapshots in a snapshot matrix, which we then use to generate the ROM basis functions (modes), as described in Section~\ref{sec:pod}. 
The snapshot average is used as the first mode, which satisfies the boundary conditions. 
We subtract the first mode from the snapshots, and solve an eigenvalue problem to find the dominant modes of these adjusted snapshots. 
For all the ROMs investigated in this section (i.e., G-ROM, DF-ROM, and HODF-ROM), we consider $r=4$.
These modes are used to generate the ROM-mass, ROM-stiffness, and ROM-nonlinear matrices, as well as the other vectors required to compute the lift and drag coefficients in the ROM online phase. 

\subsubsection{Online Phase} 
\label{sec:numerical-results-roms-online}

In the online phase of the EFR-SCM-ROM algorithm, we run the DF-ROM and HODF-ROM-$m$ that were built in Section~\ref{sec:numerical-results-roms-offline}. 
To avoid any significant temporal error in the ROM discretization, in all our tests we choose the same timestep as that used in the DNS, i.e., $\Delta t=0.002$. 
As discussed in Section~\ref{sec:efr-scm-rom}, we use the same number of collocation points and corresponding weights as those used in the DNS phase. 
For each sample point, the ROM initial condition at $t=7$ is the $L^2$ projection  into the ROM space of the corresponding FE solution at $t=7$. 
We also use the backward Euler method to construct the ROM initial condition at $t=7.002$. 
For each collocation point, we apply the DF-ROM and HODF-ROM-$m$ algorithms described in Section \ref{sec:efr-scm-rom}. 
In the EFR-SCM-ROM algorithm, we use the scaling $\chi = \Delta t$, which is commonly used in the FE discretization of EFR models~\cite{ervin2012numerical}.
We start the ROMs from the initial time $t=7.002$ (now called $t=0$) and we run them up to $t=10$. 
Since the ROM initial condition was taken at $t=7$, we plot all the ROM results on the time interval $t \in [7,17]$.
For each collocation point, we compute the energy, lift, and drag as our ROM quantities of interest. 

\subsubsection{Computational Efficiency}
\label{sec:numerical-results-roms-computational-efficiency}

As described in Section~\ref{sec:numerical-results-dns}, to generate the DNS data we integrate in time a system of equations with $35020$ velocity and $25974$ pressure degrees of freedom.
We need to integrate this high-dimensional system of equations for each collocation point: $N_{sc}=65$ for the constant 1D random viscosity case (see Section~\ref{sec:numerical-results-1d}) and $N_{sc}=801$ for the variable 5D random viscosity case (see Section~\ref{sec:numerical-results-5d}).

As described in Section~\ref{sec:numerical-results-roms}, to generate the DF-ROM and HODF-ROM-$m$ (i.e., the ROM offline phase in which the DF-ROM and HODF-ROM are built), we integrate in time the same high-dimensional system as that used to generate the DNS data.
We emphasize, however, that in this ROM offline phase we use significantly fewer collocation points: $N_{train}=9$ (instead of $N_{sc}=65$) for the constant 1D random viscosity case (see Section~\ref{sec:numerical-results-1d}) and $N_{train}=11$  (instead of $N_{sc}=801$) for the variable 5D random viscosity case (see Section~\ref{sec:numerical-results-5d}). 
Furthermore, in the ROM online phase, we integrate in time a system of equations with only $4$ degrees of freedom.

In conclusion, both the DF-ROM and the HODF-ROM provide considerable computational savings compared to the brute force DNS.

\subsection{Constant 1D Random Viscosity}
\label{sec:numerical-results-1d}

In this section, we consider the stochastic incompressible non-stationary NSE~\eqref{eq:NSE} subject to the constant random viscosity 
$$\nu_i=\nu_0\left(1+\frac{y_i}{10}\right) ,$$
where $\nu_0=8\times 10^{-4}$ and $y_i\in[-1, 1]$ is the $i^{th}$ outcome of a uniform random variable $y$ with mean zero and variance $\frac12$. 

\begin{figure}[h!]
	\begin{center}
		\includegraphics[width=0.7\textwidth,height=0.22\textwidth, trim=0 0 0 0, clip]{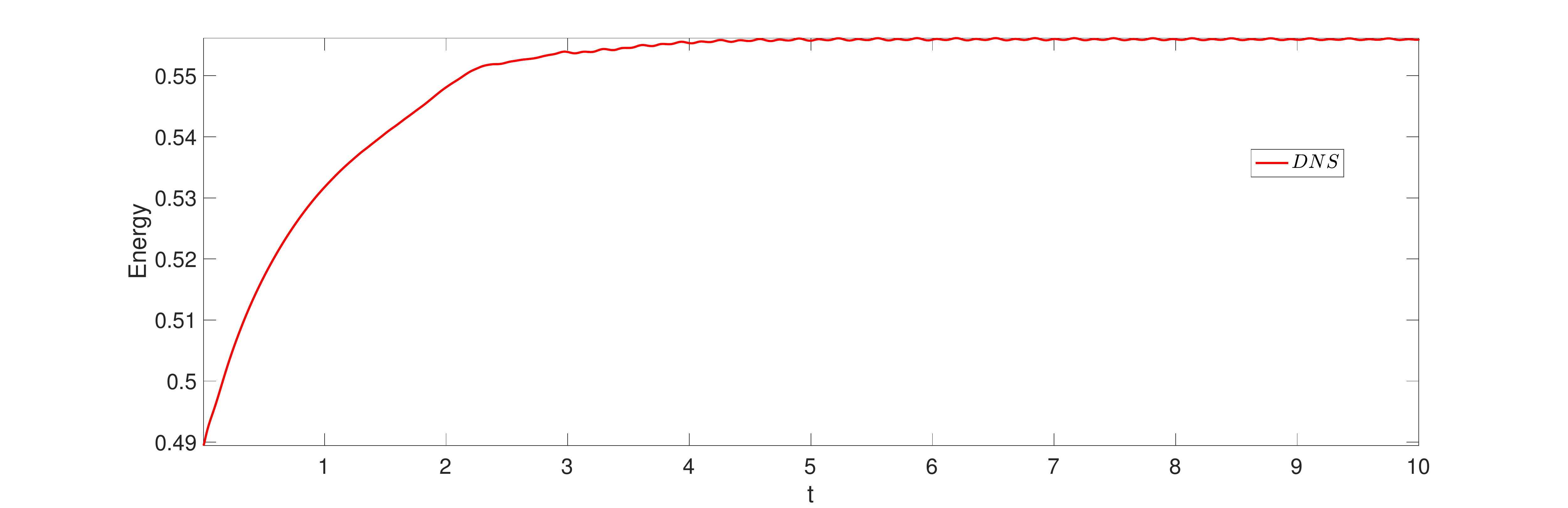}
	\end{center}
	\caption{
		Constant 1D random viscosity: 
		Plot of DNS energy coefficient vs. time.
		\label{cons_vis_dns_energy}
	}
\end{figure}

\begin{figure}[h!]
	\begin{center}
		\includegraphics[width=0.7\textwidth,height=0.22\textwidth, trim=0 0 0 0, clip]{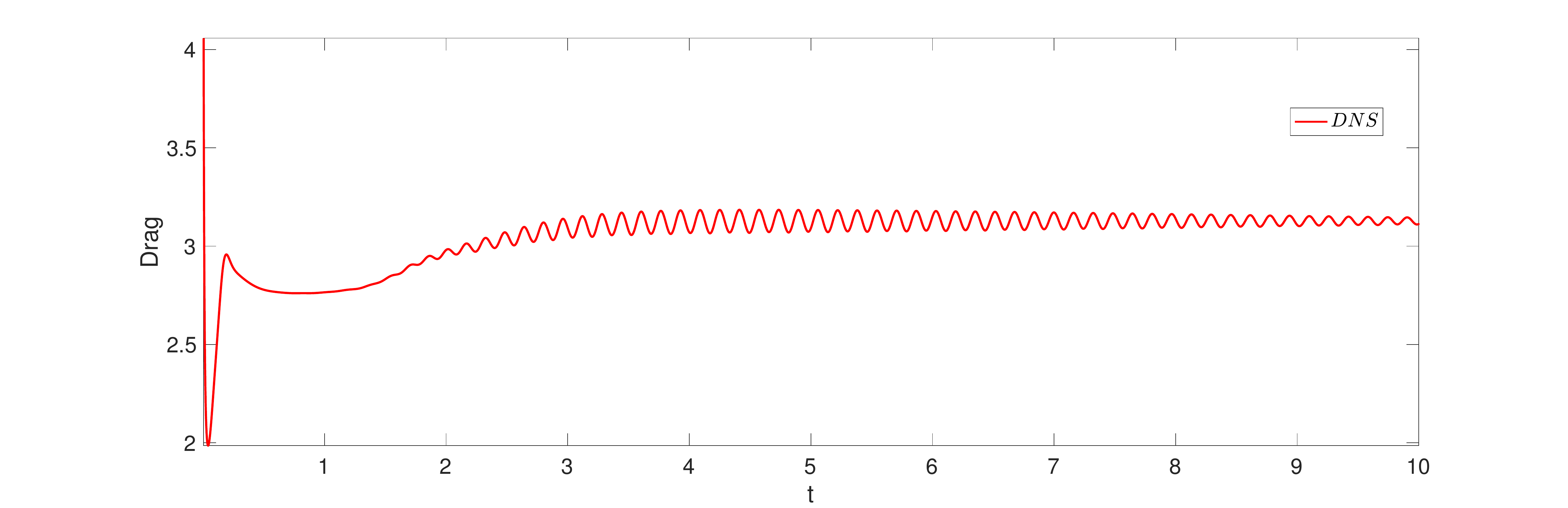}
	\end{center}
	\caption{
		Constant 1D random viscosity: 
		Plot of DNS drag coefficient vs. time.
		\label{cons_vis_dns_drag}
	}
\end{figure}

\begin{figure}[h!]
	\begin{center}
		\includegraphics[width=0.7\textwidth,height=0.22\textwidth, trim=0 0 0 0, clip]{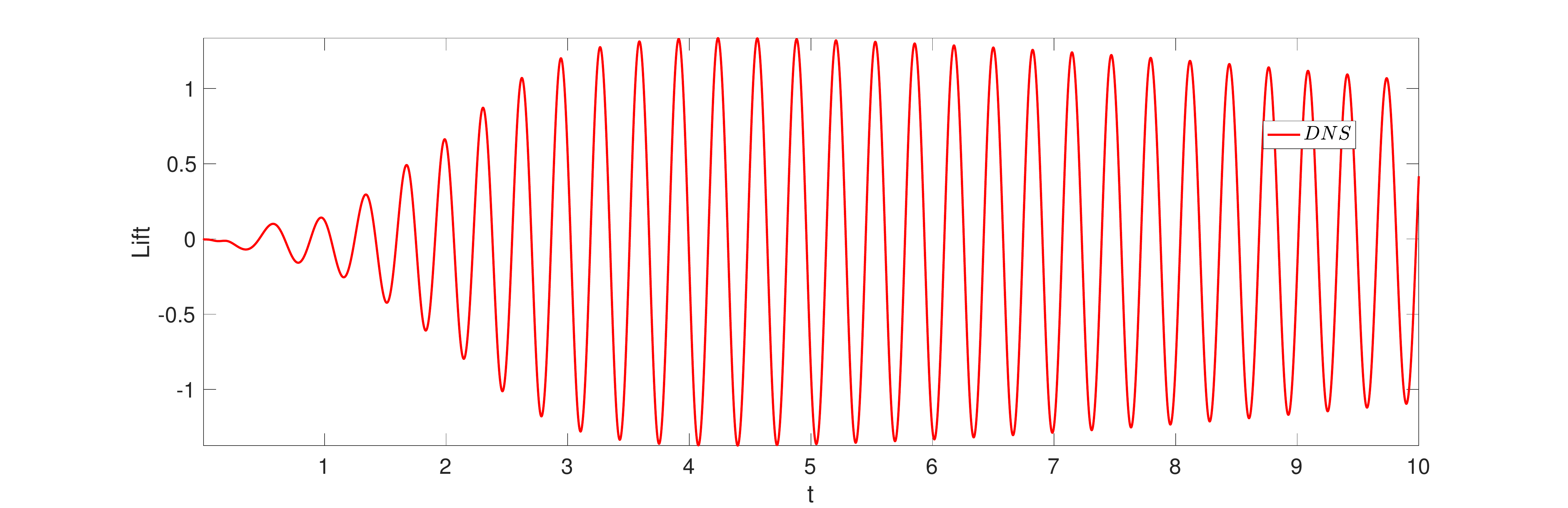}
	\end{center}
	\caption{
		Constant 1D random viscosity: 
		Plot of DNS lift coefficient vs. time.
		\label{cons_vis_dns_lift}
	}
\end{figure}

To generate the DNS data, we use $N_{sc}=65$ collocation points and their corresponding weights. 
We run the DNS for each of the collocation points and we compute the energy, lift, and drag, which are shown in Figures \ref{cons_vis_dns_energy}-\ref{cons_vis_dns_lift} as our DNS quantities of interest.

To generate the ROM data, we use $N_{train}=9$ collocation points and their corresponding weights. 
We run the ROMs for each of the collocation points and we compute the energy coefficients.

\begin{figure}[h!]
	\begin{center}
		\includegraphics[width=0.7\textwidth,height=0.22\textwidth, trim=0 0 0 0, clip]{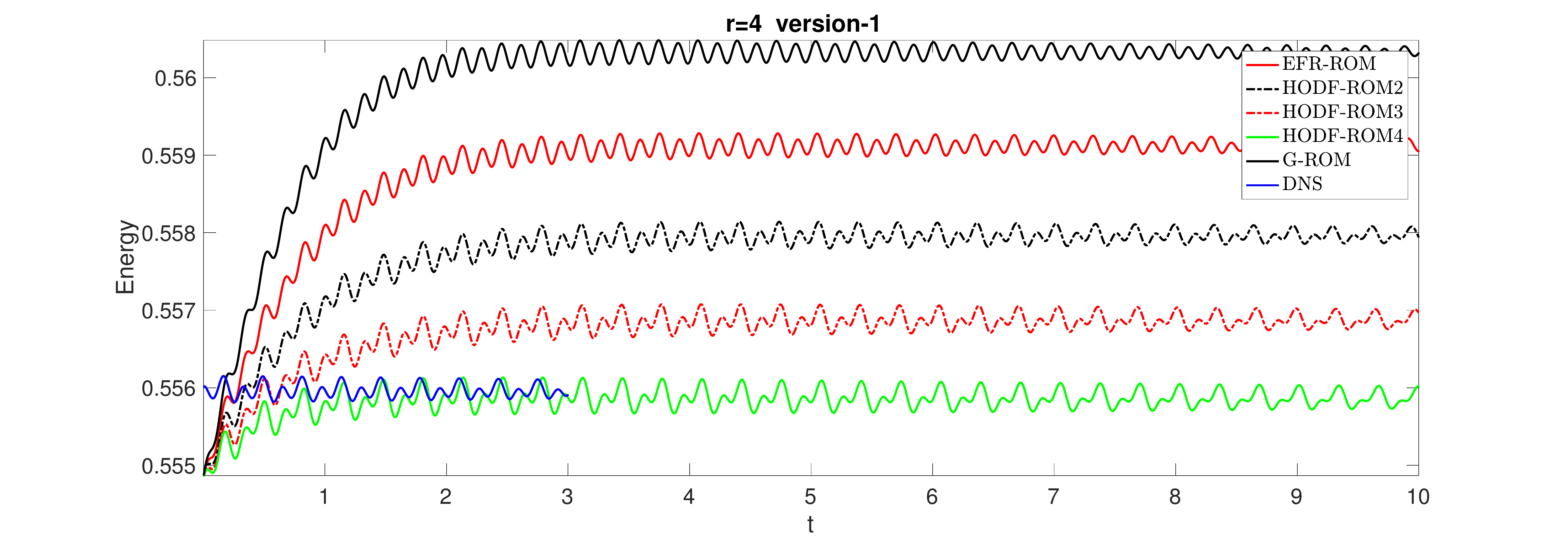}		
	\end{center}
	\caption{
		Constant 1D random viscosity: 
		Plot of DNS energy coefficient vs. time for $\delta=7.5\times 10^{-3}$.
		\label{fig:1d-v1-1}
	}
\end{figure}

\begin{figure}[h!]
	\begin{center}
		\includegraphics[width=0.7\textwidth,height=0.22\textwidth, trim=0 0 0 0, clip]{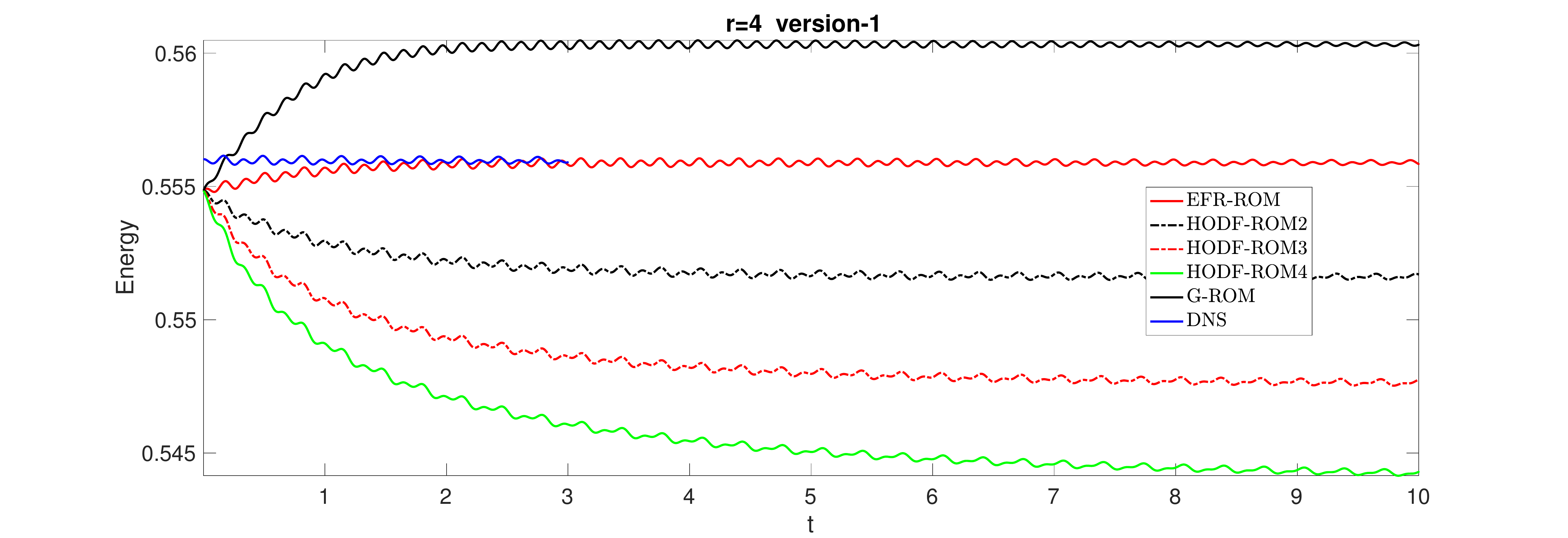}
	\end{center}
	\caption{
		Constant 1D random viscosity: 
		Plot of DNS energy coefficient vs. time for $\delta=1\times 10^{-2}$.
		\label{fig:1d-v1-2}
	}
\end{figure}

First, we consider version 1 of the HODF, i.e., equation~\eqref{eqn:hodf-1}.
If Figures~\ref{fig:1d-v1-1}--\ref{fig:1d-v1-2}, we plot the time evolution of the energy coefficients of the four models investigated in this section: DNS, G-ROM, DF-ROM, and version 1 of HODF-ROM-$m$ with $m=2, 3, 4$.
We note that the time evolution of the drag and lift coefficients of all the models follow the same trends as  the time evolution of the energy coefficients; for clarity of presentation, we do not include the former.
For the DF-ROM, and version 1 of HODF-ROM-$m$, we consider the following $\delta$ values: $\delta = 5 \times 10^{-3}, 6 \times 10^{-3}, 7 \times 10^{-3}, 7.5 \times 10^{-3}, 8 \times 10^{-3}, 1 \times 10^{-2}, 2 \times 10^{-2}$.
However, for clarity of presentation, we include results only for $\delta = 7.5 \times 10^{-3}, 1 \times 10^{-2}$.

The plots in Figures~\ref{fig:1d-v1-1}--\ref{fig:1d-v1-2} show that by varying $\delta$ one can get optimal accuracy for both the DF-ROM and HODF-ROM-$m$.  
Figure~\ref{fig:1d-v1-1} shows that, for small $\delta$ values (i.e., $\delta = 7.5 \times 10^{-3}$), the HODF-ROM-$m$ with $m=4$ performs dramatically better than G-ROM and significantly better than DF-ROM.
Furthermore, the $m=4$ value yields the best results for the HODF-ROM-$m$.
Figure~\ref{fig:1d-v1-2} shows that, for large $\delta$ values (i.e., $\delta=1 \times 10^{-2}$), the situation is reversed: DF-ROM is more accurate than HODF-ROM-$m$.
We emphasize, however, that even in this case both DF-ROM and HODF-ROM-$m$ are significantly more accurate than G-ROM.

Overall, the plots in Figures~\ref{fig:1d-v1-1}--\ref{fig:1d-v1-2} yield the following conclusions:
For small $\delta$ values, HODF-ROM-$m$ (with high $m$ values) is the most accurate.
For large $\delta$ values, DF-ROM is the most accurate.
For a fixed $\delta$ value, the tuning of the extra parameter $m$ in HODF-ROM-$m$ allows it to perform better than DF-ROM. 

\bigskip

\begin{figure}[h!]
	\begin{center}
		\includegraphics[width=0.7\textwidth,height=0.22\textwidth, trim=0 0 0 0, clip]{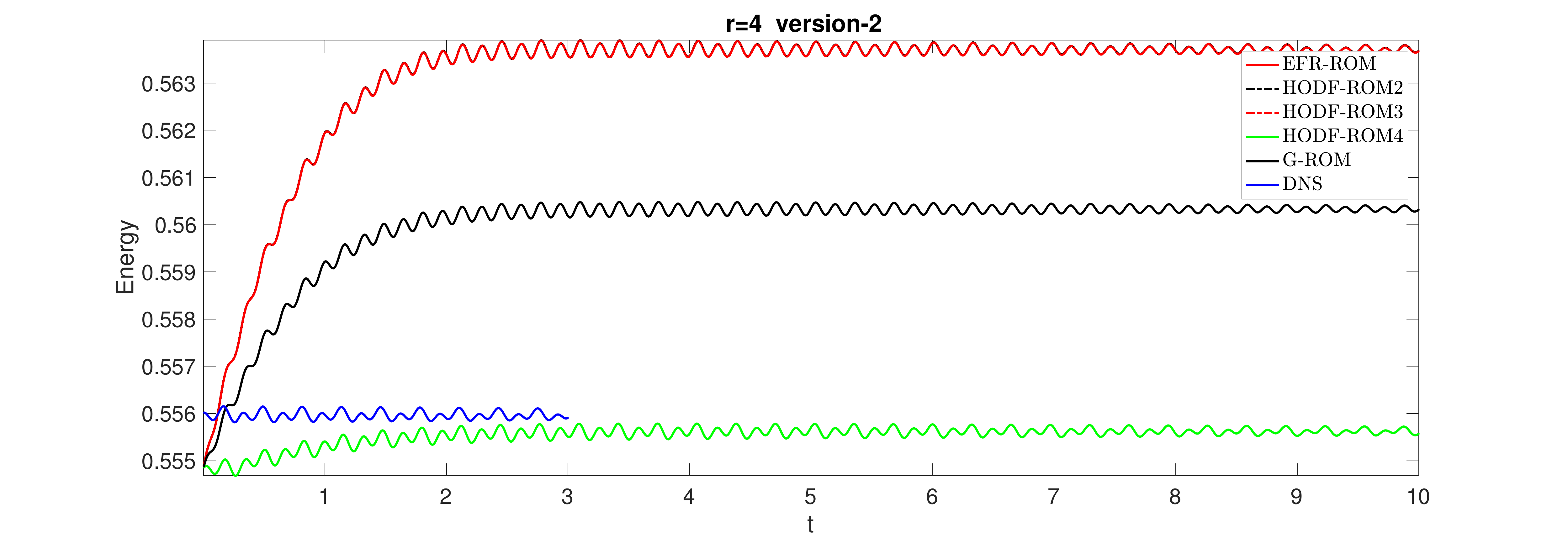}
	\end{center}
	\caption{
		Constant 1D random viscosity: 
		Plot of DNS energy coefficient vs. time for $\delta=1.5\times 10^{-7}$.
		\label{fig:1d-v2-1}
	}
\end{figure}

\begin{figure}[h!]
	\begin{center}
		\includegraphics[width=0.7\textwidth,height=0.22\textwidth, trim=0 0 0 0, clip]{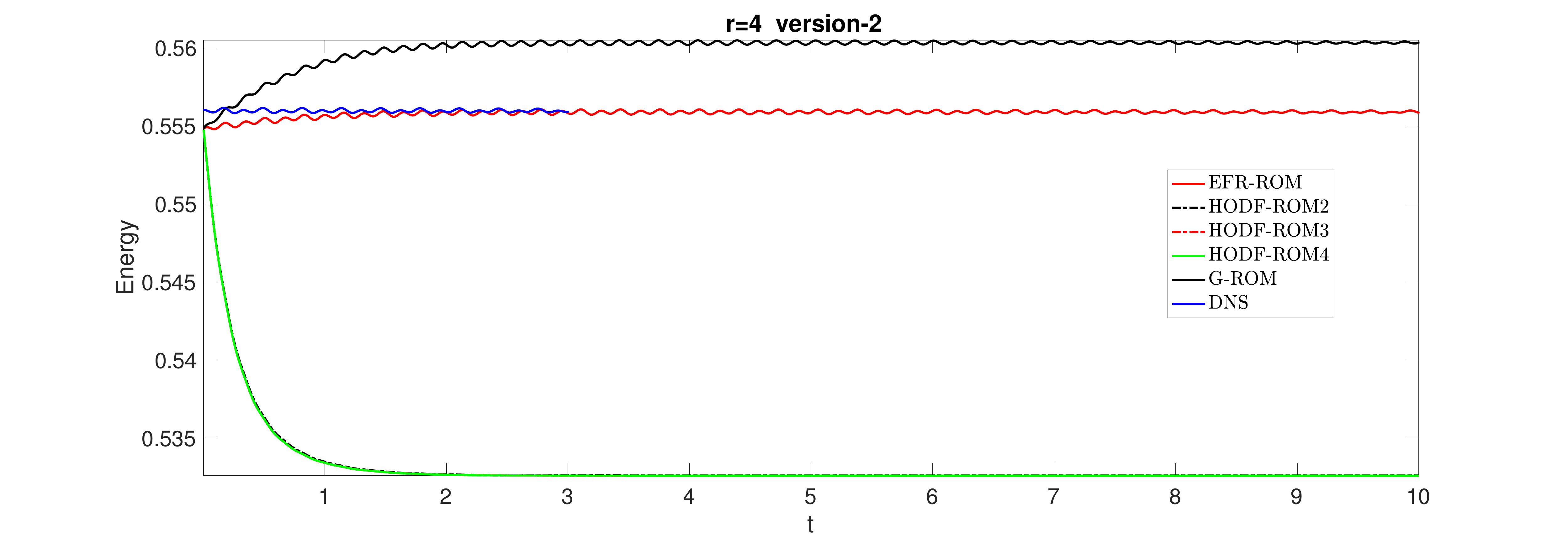}
	\end{center}
	\caption{
		Constant 1D random viscosity: 
		Plot of DNS energy coefficient vs. time for $\delta=1\times 10^{-2}$.
		\label{fig:1d-v2-2}
	}
\end{figure}

Next, we consider version 2 of the HODF, i.e., equation~\eqref{eqn:hodf-2}.
If Figures~\ref{fig:1d-v2-1}--\ref{fig:1d-v2-2}, we plot the time evolution of the energy coefficients of the four models investigated in this section: DNS, G-ROM, DF-ROM, and version 2 of HODF-ROM-$m$ with $m=2, 3$, and $4$.
For the DF-ROM and version 2 of HODF-ROM-$m$, we consider the following $\delta$ values: $\delta = 1 \times 10^{-7}, 1.5 \times 10^{-7}, 2 \times 10^{-7}, 1 \times 10^{-6}, 1 \times 10^{-5}, 1 \times 10^{-3}, 5 \times 10^{-3}, 1 \times 10^{-2}, 2 \times 10^{-2}$.
However, for clarity of presentation, we include results only for $\delta = 1.5 \times 10^{-7}, 1 \times 10^{-2}$.

The plots in Figures~\ref{fig:1d-v2-1}--\ref{fig:1d-v2-2} (for version 2 of HODF-ROM-$m$) yield the same conclusions as the plots in Figures~\ref{fig:1d-v1-1}--\ref{fig:1d-v1-2} (for version 1 of HODF-ROM-$m$):
For small $\delta$ values, HODF-ROM-$m$ (with high $m$ values) is the most accurate.
For large $\delta$ values, DF-ROM is the most accurate.
For a fixed $\delta$ value, the tuning of the extra parameter $m$ in HODF-ROM-$m$ allows it to perform better than DF-ROM.


\subsection{Variable 5D Random Viscosity}
\label{sec:numerical-results-5d}

In this section, we consider the time-dependent incompressible stochastic NSE~\eqref{eq:NSE} with a random viscosity $\nu({x},{y})$, where ${y}=(y_1,y_2,\cdots,y_d)\in\Gamma\subset\mathbb{R}^d$ is a higher-dimensional random variable, $\mathbb{E}[\nu](x)=\frac{c}{1000}$ for a suitable $c>0$, $\mathbb{C}ov[\nu]({x},{x^{'}})=\frac{1}{1000^2}exp\left(-\frac{({x}-{x^{'}})^2}{l^2}\right)$, and $l$ is the correlation length. 
This random field can be represented by the Karhunen-Lo\'eve expansion\\
\begin{align}
\tiny\nu({x}, {y})=\frac{1}{1000}\bigg(c+\left(\frac{\sqrt{\pi}l}{2}\right)^{\frac12}y_{1}(\omega)+\sum_{j=1}^{q}\sqrt{\xi_j}&\bigg(\sin\left(\frac{j\pi x_1}{2.2}\right)\sin\left(\frac{j\pi x_2}{0.41}\right)y_{2j}(\omega)\nonumber\\&+\cos\left(\frac{j\pi x_1}{2.2}\right)\cos\left(\frac{j\pi x_2}{0.41}\right)y_{2j+1}(\omega)\bigg), \label{eq:var-vis}
\end{align}
in which the infinite series is truncated up to the first $q$ terms. 
The uncorrelated random variables $y_{j}$ have zero mean and unit variance, and the eigenvalues are equal to
$$\sqrt{\xi_j}=(\sqrt{\pi}l)^{\frac12}exp\left(-\frac{(j\pi l)^2}{8}\right).$$
For our test problem, we consider the random variables $y_{j}(\omega)\in[-\sqrt{3},\sqrt{3}]$, the correlation length $l=0.01$, $N=5$, $c=1$, and $q=2$.

To generate the DNS data, we use $N_{sc}=801$ collocation points and their corresponding weights.

\begin{figure}[h!]
	\begin{center}
		\includegraphics[width=0.7\textwidth,height=0.22\textwidth, trim=0 0 0 0, clip]{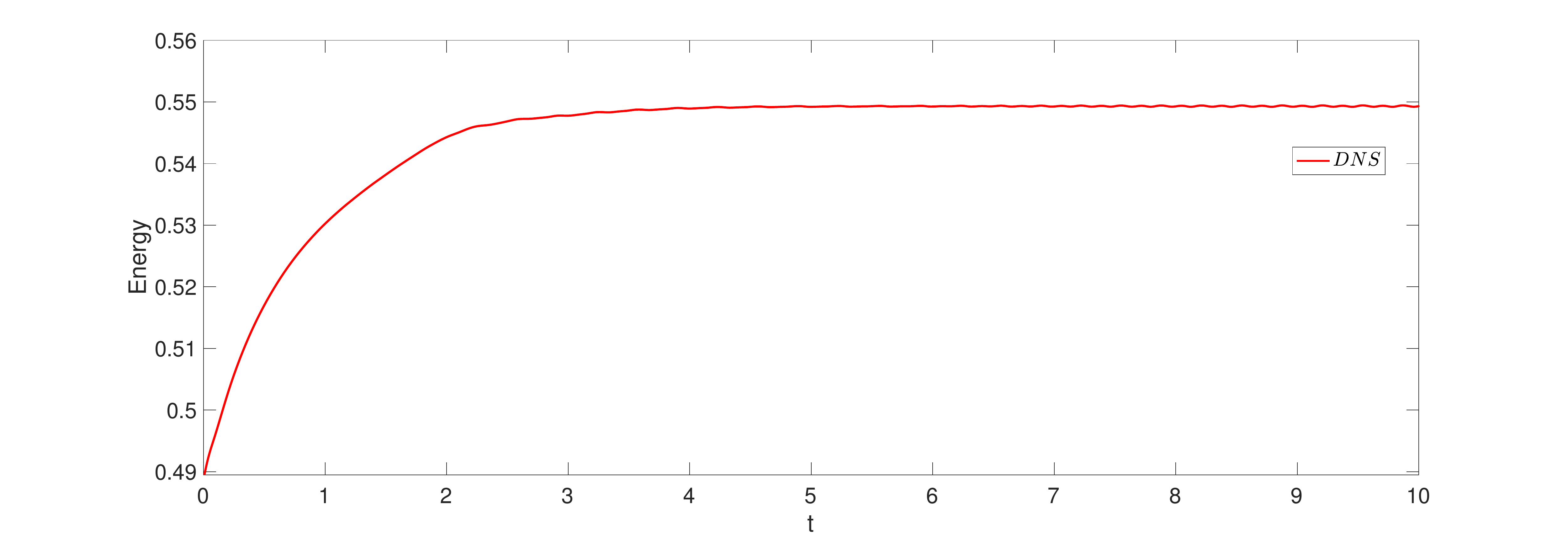}
	\end{center}
	\caption{
		Variable 5D random viscosity: 
		Plot of DNS energy coefficient vs. time.
		\label{bn1}
	}
\end{figure}

\begin{figure}[h!]
	\begin{center}
		\includegraphics[width=0.7\textwidth,height=0.22\textwidth, trim=0 0 0 0, clip]{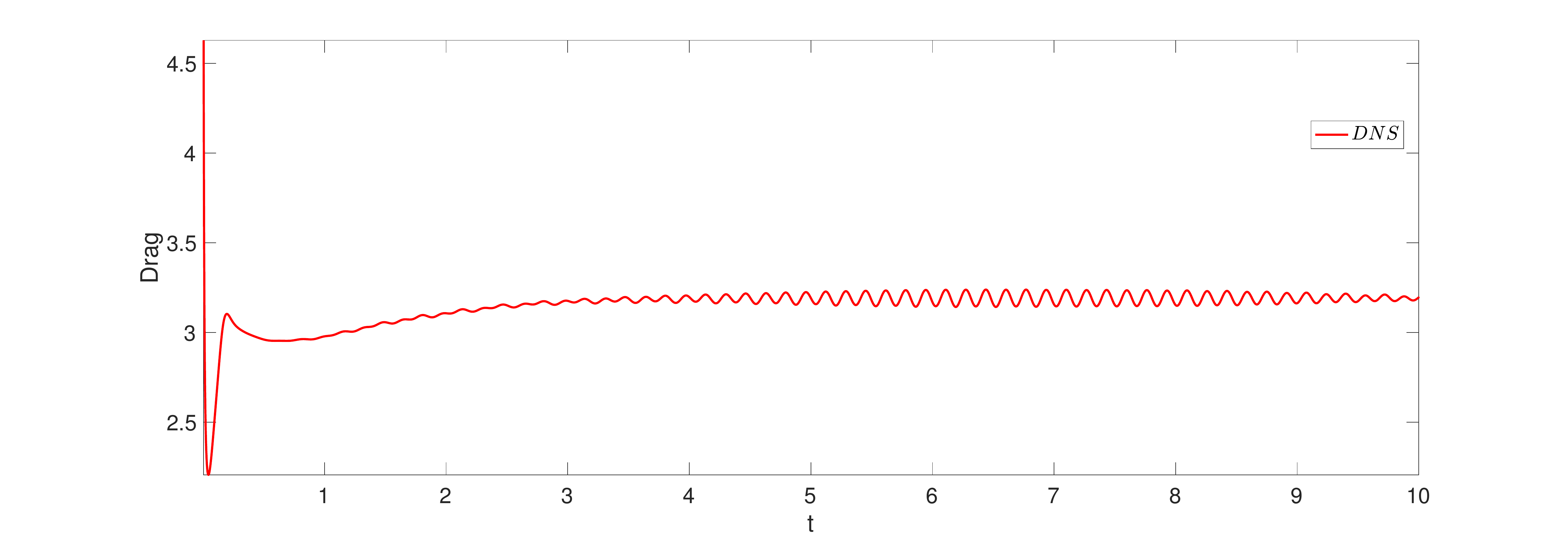}
	\end{center}
	\caption{
		Variable 5D random viscosity: 
		Plot of DNS drag coefficient vs. time.
		\label{bn2}
	}
\end{figure}

\begin{figure}[h!]
	\begin{center}
		\includegraphics[width=0.7\textwidth,height=0.22\textwidth, trim=0 0 0 0, clip]{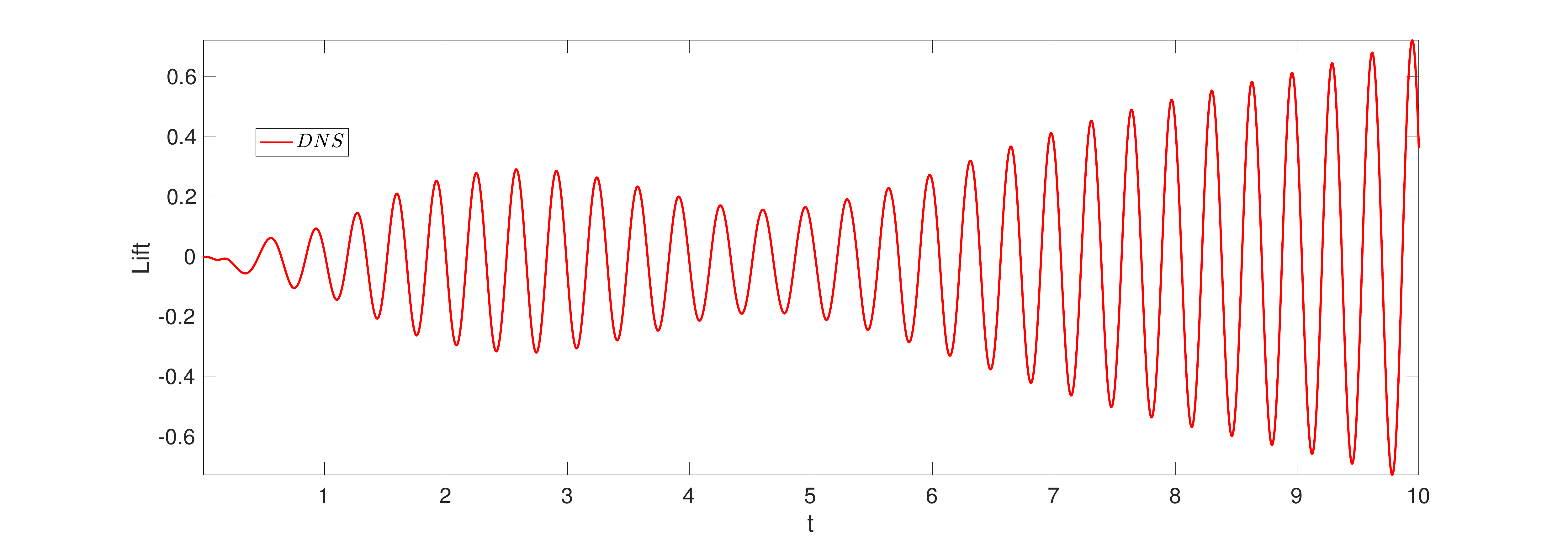}
	\end{center}
	\caption{
		Variable 5D random viscosity: 
		Plot of DNS lift coefficient vs. time.
		\label{bn6}
	}
\end{figure}

We run the DNS for each of the collocation points and we compute the energy, lift, and drag (which are shown in Figures~\ref{bn1}-\ref{bn6}) as our DNS quantities of interest.

To generate the ROM data, we use $N_{train}=11$ collocation points.

We run the ROMs for each of the collocation points and we compute the energy coefficients.

\begin{figure}[h!]
	\begin{center}
		\includegraphics[width=0.7\textwidth,height=0.22\textwidth, trim=0 0 0 0, clip]{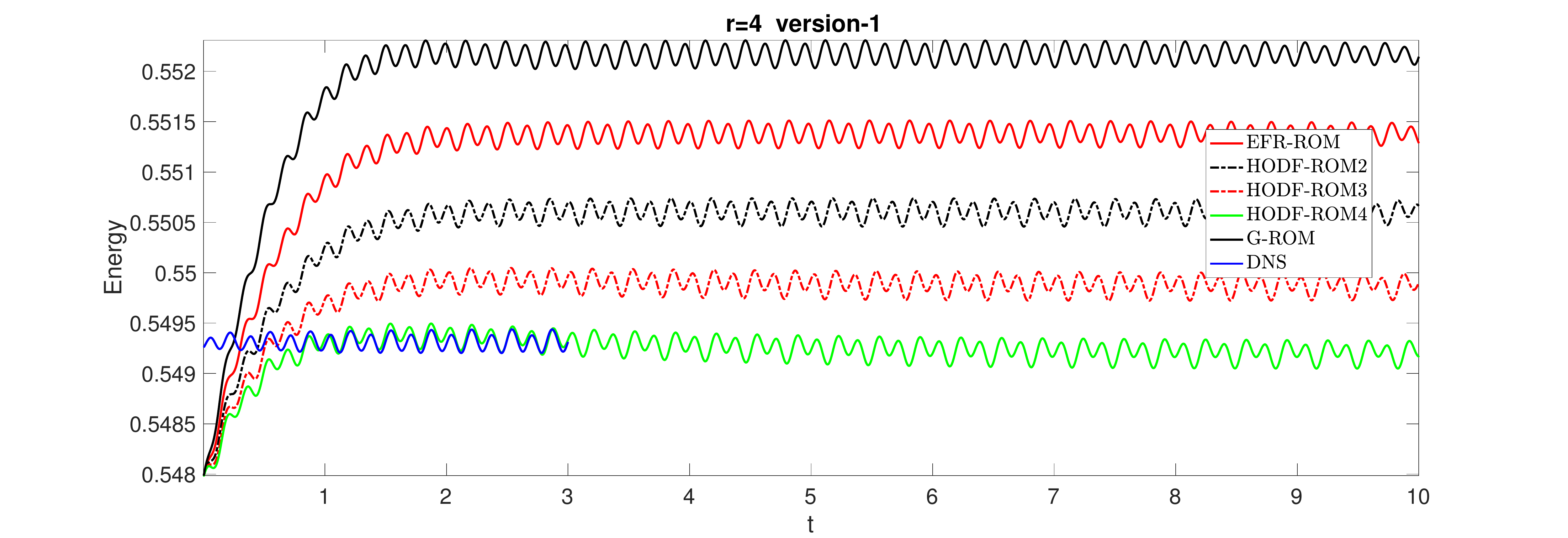}
	\end{center}
	\caption{
		Variable 5D random viscosity: 
		Plot of DNS energy coefficient vs. time for $\delta=9.6\times 10^{-3}$.
		\label{fig:5d-v1-1}
	}
\end{figure}

\begin{figure}[h!]
	\begin{center}
		\includegraphics[width=0.7\textwidth,height=0.22\textwidth, trim=0 0 0 0, clip]{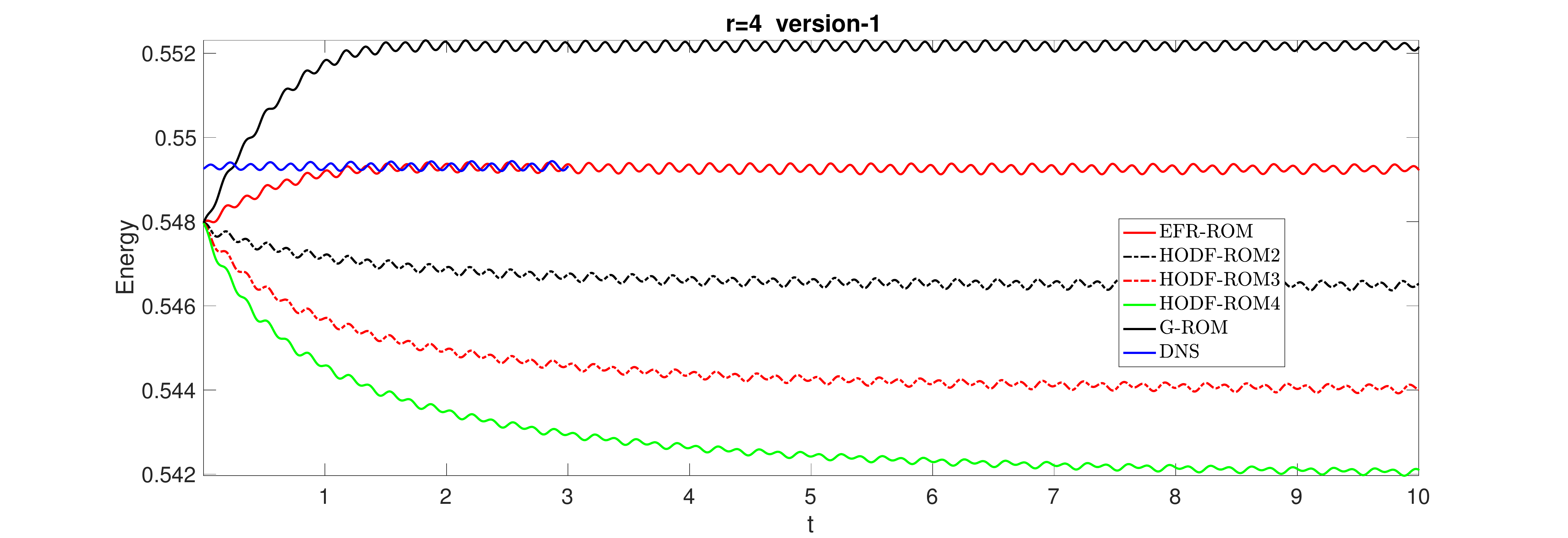}
	\end{center}
	\caption{
		Variable 5D random viscosity: 
		Plot of DNS energy coefficient vs. time for $\delta=1.18\times 10^{-2}$.
		\label{fig:5d-v1-2}
	}
\end{figure}

First, we consider version 1 of the HODF, i.e., equation~\eqref{eqn:hodf-1}.
If Figures~\ref{fig:5d-v1-1}--\ref{fig:5d-v1-2}, we plot the time evolution of the energy coefficients of DNS, G-ROM, DF-ROM, and version 1 of HODF-ROM-$m$ with $m=2, 3$, and $4$.
As in the constant viscosity 1D case, the time evolution of the drag and lift coefficients of all the models follow the same trends as  the time evolution of the energy coefficients; for clarity of presentation, we do not include the former.
For the DF-ROM, and version 1 of HODF-ROM-$m$, we consider the following $\delta$ values: $\delta = 1 \times 10^{-3}, 5 \times 10^{-3}, 6 \times 10^{-3}, 7 \times 10^{-3}, 8 \times 10^{-3}, 9 \times 10^{-3}, 9.5 \times 10^{-3}, 9.6 \times 10^{-3}, 1 \times 10^{-2}, 1.18 \times 10^{-2}, 2 \times 10^{-2}$.
However, for clarity of presentation, we include results only for $\delta = 9.6 \times 10^{-3}, 1.18 \times 10^{-2}$.

As in the constant viscosity 1D case, the plots in Figures~\ref{fig:5d-v1-1}--\ref{fig:5d-v1-2} show that by varying $\delta$ one can get optimal accuracy for both the DF-ROM and HODF-ROM-$m$.  
Figure~\ref{fig:5d-v1-1} shows that, for small $\delta$ values (i.e., $\delta = 9.6 \times 10^{-3}$), the HODF-ROM-$m$ with $m=4$ performs dramatically better than G-ROM and significantly better than DF-ROM.
Furthermore, the $m=4$ value yields the best results for the HODF-ROM-$m$.
Figure~\ref{fig:5d-v1-2} shows that, for large $\delta$ values (i.e., $\delta=1.18 \times 10^{-2}$), the situation is reversed: DF-ROM is more accurate than HODF-ROM-$m$.
We emphasize, however, that even in this case both DF-ROM and HODF-ROM-$m$ are significantly more accurate than G-ROM.

Overall, the plots in Figures~\ref{fig:5d-v1-1}--\ref{fig:5d-v1-2} yield the same conclusions as in the constant viscosity 1D case:
For small $\delta$ values, HODF-ROM-$m$ (with high $m$ values) is the most accurate.
For large $\delta$ values, DF-ROM is the most accurate.
For a fixed $\delta$ value, the tuning of the extra parameter $m$ in HODF-ROM-$m$ allows it to perform better than DF-ROM. 

\bigskip

\begin{figure}[h!]
	\begin{center}
		\includegraphics[width=0.7\textwidth,height=0.22\textwidth, trim=0 0 0 0, clip]{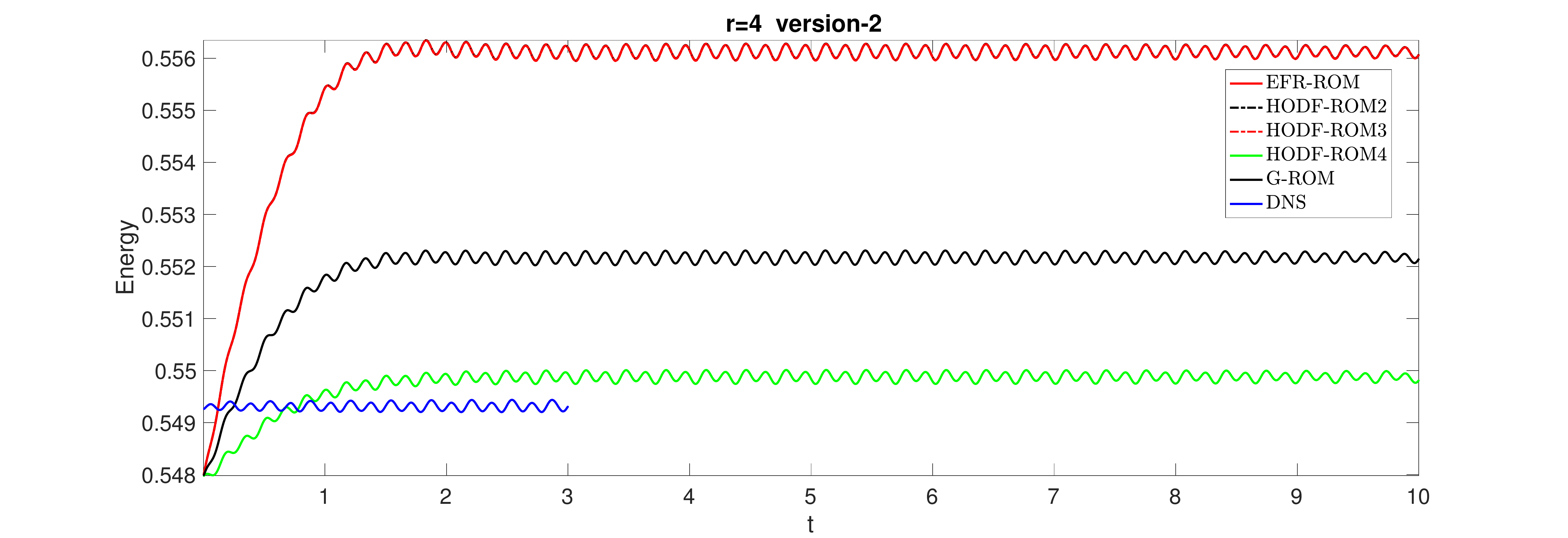}
	\end{center}
	\caption{
		Variable 5D random viscosity: 
		Plot of DNS energy coefficient vs. time for $\delta=2\times 10^{-7}$.
		\label{fig:5d-v2-1}
	}
\end{figure}

\begin{figure}[h!]
	\begin{center}
		\includegraphics[width=0.7\textwidth,height=0.22\textwidth, trim=0 0 0 0, clip]{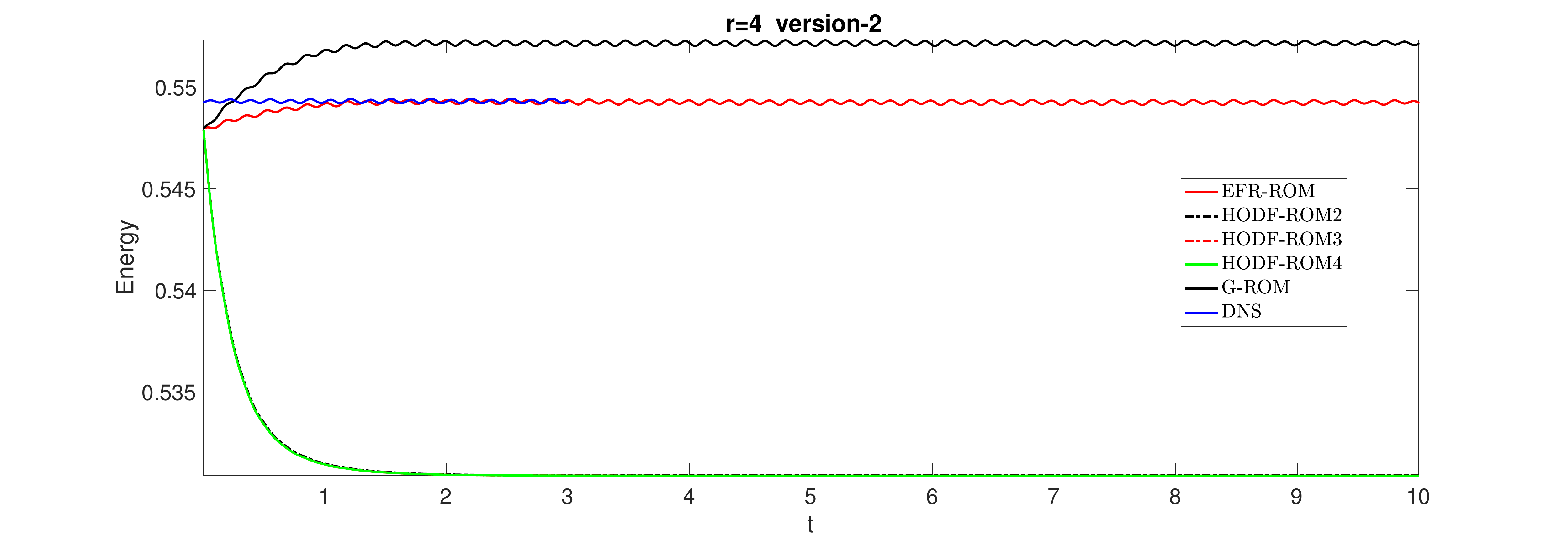}
	\end{center}
	\caption{
		Variable 5D random viscosity: 
		Plot of DNS energy coefficient vs. time for $\delta=1.18\times 10^{-2}$.
		\label{fig:5d-v2-2}
	}
\end{figure}

Next, we consider version 2 of the HODF, i.e., equation~\eqref{eqn:hodf-2}.
If Figures~\ref{fig:5d-v2-1}--\ref{fig:5d-v2-2}, we plot the time evolution of the energy coefficients of DNS, G-ROM, DF-ROM, and version 2 of HODF-ROM-$m$ with $m=2, 3$, and $4$.
For the DF-ROM and version 2 of HODF-ROM-$m$, we consider the following $\delta$ values: $\delta = 1 \times 10^{-7}, 2 \times 10^{-7}, 1 \times 10^{-6}, 1 \times 10^{-5}, 5 \times 10^{-3}, 1 \times 10^{-2}, 1.18 \times 10^{-2}, 2 \times 10^{-2}$.
However, for clarity of presentation, we include results only for $\delta = 2 \times 10^{-7}, 1.18 \times 10^{-2}$.

The plots in Figures~\ref{fig:5d-v2-1}--\ref{fig:5d-v2-2} (for version 2 of HODF-ROM-$m$) yield the same conclusions as the plots in Figures~\ref{fig:5d-v1-1}--\ref{fig:5d-v1-2} (for version 1 of HODF-ROM-$m$):
For small $\delta$ values, HODF-ROM-$m$ (with high $m$ values) is the most accurate.
For large $\delta$ values, DF-ROM is the most accurate.
For a fixed $\delta$ value, the tuning of the extra parameter $m$ in HODF-ROM-$m$ allows it to perform better than DF-ROM. 

\section{Conclusions}
\label{sec:conclusions}

We proposed a nonintrusive filter-based stabilization of ROMs for UQ of the time-dependent NSE in convection-dominated regimes.
For the ROM stabilization, we put forth a new evolve-filter-relax (EFR) algorithm to attenuate the numerical oscillations that standard ROMs generally display in convection-dominated flows.
Furthermore, we proposed a novel high-order ROM differential filter (HODF) to \textcolor{black}{perform} the spatial filtering in the evolve-filter-relax algorithm.
For the UQ component, we used the SCM.
We emphasize that the entire stabilized SCM-ROM framework that we proposed is nonintrusive and can be easily used in conjunction with legacy flow solvers.
We also note that, to our knowledge, the SCM-ROM for time-dependent NSE, the EFR algorithm for ROMs, and the HODF are new.

We performed a numerical investigation of the new SCM-EFR-ROM equipped with the DFDF~\eqref{eqn:df} and the HODF~\eqref{eqn:hodf-1} with $m=2, 3$, and $4$.
As a test problem, we considered the 2D flow past a circular cylinder with a random viscosity that yields a \textcolor{black}{random} Reynolds number with mean $Re=100$.
We considered two cases for the random viscosity: a 1D constant viscosity case and a 5D variable viscosity case. 
Our numerical investigation yielded the following general conclusions:
For small $\delta$ values, HODF-ROM-$m$ (with high $m$ values) was the most accurate.
For large $\delta$ values, DF-ROM was the most accurate.
For a fixed $\delta$ value, the tuning of the extra parameter $m$ in HODF-ROM-$m$ allowed it to perform better than DF-ROM. 

\medskip

These first steps in the numerical investigation of the new SCM-EFR-ROM framework are encouraging.
There are, however, several research avenues that we plan to pursue for a better understanding of the capabilities and limitations of the SCM-EFR-ROM.

First, we plan to continue the comparison of the DF vs. HODF comparison, both in a deterministic and a stochastic setting.
In particular, we will investigate the effect of the two ROM spatial filters on the ROM spatial structures.
To this end, in addition to the time evolution of the energy coefficients, we will consider \textcolor{black}{alternative} criteria in our numerical investigation, such as the time evolution of the ROM coefficients, $a_{i}$.

We also plan to investigate the new EFR-SCM-ROM framework in more challenging flows, such as the NSE at higher Reynolds numbers.
In particular, to study the limitations of the proposed EFR-SCM-ROMs, we will gradually increase the Reynolds number until the accuracy of the ROM approximation exceeds the prescribed tolerance.
We also plan to perform a sensitivity study of the ROM spatial filter radius $\delta$ that we used in the DF and HODF with respect to the Reynolds number.
{Lastly, we plan to examine the effect of the spatial filters on the accuracy of the SCM for long time evolutions. It is known that as the final time grows, the error from the SCM will also grow~\cite{tran2014convergence, WK06}. Whether or not the spatial filters mitigate this effect is a subject of future work.}

\clearpage

\bibliographystyle{plain}
\bibliography{traian}
\end{document}